\documentclass[11pt, reqno]{amsart}
\usepackage{amssymb}
\usepackage{lmodern}
\usepackage{graphicx}
\usepackage{physics}
\usepackage{amssymb}
 \usepackage{amsthm}
\usepackage{epsf}
\usepackage{amsbsy,amsmath}
\usepackage{mathtools}
\usepackage{upgreek}
\usepackage{mathrsfs}
\usepackage{amsfonts}
\usepackage{amssymb}
\usepackage{float}
\usepackage{enumitem}
\usepackage{eucal}
\usepackage{graphics,mathrsfs}
\usepackage{amsthm}
\usepackage{secdot}
\usepackage{esint}
\usepackage{hyperref}
\usepackage{varwidth}
\usepackage{tasks}
\usepackage{abstract}
\theoremstyle{plain}
\newtheorem{theorem}{Theorem}[section]
\newtheorem{cor}{Corollary}[section]

\theoremstyle{definition}
\newtheorem{definition}{Definition}[section] % This numbers definitions based on sections
\newtheorem{lemma}{Lemma}[section]
\newtheorem{remark}{Remark}[section]

\allowdisplaybreaks

\usepackage{xcolor}
\usepackage{bbm}
\usepackage{cite}
\usepackage{comment}

\def\div{{\operatorname{div}}}

\renewcommand{\d}{\:\! \mathrm{d}}

\long\def\symbolfootnote[#1]#2{\begingroup
\def\thefootnote{\fnsymbol{footnote}}\footnote[#1]{#2}\endgroup}
\numberwithin{equation}{section}
\begin{document}
\title[Least energy sign-changing solutions to Kirchhoff-type problems] {Least energy sign-changing solutions to local-nonlocal Kirchhoff-type problems}
\author{Souvik Bhowmick}
\address[Souvik Bhowmick]{Department of Mathematics, National Institute of Technology Calicut, Kozhikode, Kerala, India - 673601}
\email{souvikbhowmick2912@gmail.com / souvik\_p230197ma@nitc.ac.in}
\author{Sekhar Ghosh}
\address[Sekhar Ghosh]{Department of Mathematics, National Institute of Technology Calicut, Kozhikode, Kerala, India - 673601}
\email{sekharghosh1234@gmail.com / sekharghosh@nitc.ac.in}

\author{Kamel Saoudi}
\address[Kamel Saoudi]{Basic and Applied Scientific Research Center, Imam Abdulrahman Bin Faisal University, Dammam, 34212, Saudi Arabia \newline and \newline 
College of Sciences at Dammam, Imam Abdulrahman Bin Faisal University, Dammam, 34212, Saudi Arabia}
\email{kmsaoudi@iau.edu.sa / kasaoudi@gmail.com}

\maketitle
\begin{abstract}
In this paper, we study the existence of least energy sign-changing solutions and ground state solutions for a class of Kirchhoff-type problems involving mixed local-nonlocal operators under mild assumptions on the nonlinearity. We employ a suitable variational method, the Nehari manifold method, and a Non-Nehari manifold method to establish our results. Further, we show that the least energy of a sign-changing solution is strictly greater than twice the ground state energy. Moreover, we prove the existence of the least energy sign-changing solution for a class of mixed local-nonlocal problems under the same assumptions on the nonlinearity. Finally, we prove the asymptotic property of the sign-changing solutions $u_b$ for the Kirchhoff-type problem as $b\searrow 0$.
% We extended the results of W. Shuai [J. Differential Equations, 259(4), 2015], X. H. Tang and B. Cheng [J. Differential Equations, 261(4), 2016], H. Luo, X.H. Tang and Z. Gao [J. Math. Phys. 59(4), 2018] and S. Bhowmick and S. Ghosh [arXiv:2503.04361, Preprint, 2025] for mixed local-nonlocal Kirchhoff-type problems.
\end{abstract}

%\tableofcontents
\section{Introduction and main results application}\label{sec1}

The main goal of this paper is to study the following class of problems
   \begin{equation}\label{New P}
   \begin{aligned}
       \left\{a+b \bigg(\int_\Omega|\nabla u|^2 \d x+\int_{\mathbb{R}^{2N}} \frac{|u(x)-u(y)|^2}{|x-y|^{N+2s}}\d x\d y \bigg)\right\} \mathcal{L}u&=f(u)~\text{ in }\Omega, \\
       u&=0~ \text{ in } \mathbb R^N \setminus \Omega,
       \end{aligned}
   \end{equation}
where $\Omega\subset\mathbb{R}^N$, $N=1,2,3$ is a bounded domain with smooth boundary $\partial \Omega$, $s\in(0,1)$, $a,b$ are positive constants and $\mathcal{L}:=-\Delta  +(-\Delta)^s$. Here, $(-\Delta)^s$ denotes the fractional Laplace operator defined as the Cauchy principal value, which is given by
 $$(-\Delta)^s u(x) = C({N,s}) P.V \int_{\mathbb{R}^N} \frac{u(x)-u(y)}{|x-y|^{N+2s}} \, \d y,$$
where $C({N,s})$ is a normalizing constant. For simplicity, we fix $C({N,s})=1$. The function $f: \mathbb{R}\rightarrow \mathbb{R}$ satisfies the following growth conditions:
\begin{itemize}
       \item[$(f_1)$] $f\in C(\mathbb{R}, \mathbb{R})$ and $\lim_{u\rightarrow 0}\frac{f(u)}{u}=0.$ 
       \item[$(f_2)$] There exist two constants $C_0>0$ and $q$ with $4<q<2^*=\begin{cases}
           & 6 \text{ when } N=3,\\
           & \infty \text{ when } N=1,2,
       \end{cases}$ 
       such that 
   $$|f(u)|\leq C_0(1+|u|^{q-1}),~\forall~ u\in \mathbb R.$$
       \item[$(f_3)$] $\lim_{|u|\rightarrow \infty}\frac{f(u)}{u^3}=+\infty.$ 
        \item[$(f_4)$] For all $t>0$ and $u\in \mathbb{R}\setminus \{0\}$, there exists a fixed $\theta\in(0,1)$ such that 
        \begin{equation}\label{f_4 con}
            \left\{\frac{f(u)}{u^3}-\frac{f(tu)}{(tu)^3} \right\}\operatorname{sign}(1-t)+\frac{a\theta\lambda_1|1-t^2|}{(tu)^2} \geq 0,
        \end{equation}
        where $\lambda_1$ is the first eigenvalue of $(-\Delta+(-\Delta)^s)$.
\end{itemize}
In addition to the problem \eqref{New P}, we also investigate the following problems for $b=0$ under the same assumptions on $f$.
 \begin{equation}\label{New P2}
   \begin{aligned}
       a \left(-\Delta u +(-\Delta)^s  u\right)&=f(u)~\text{ in }\Omega,\\
       u&=0~ \text{ in } \mathbb R^N \setminus \Omega.
       \end{aligned}
   \end{equation}
  It is noteworthy to mention that $f(u)=|u|^{q-2}u$ is the simplest function satisfying all assumptions $(f_1)$--$(f_4)$ for the range of $q$ as in ($f_2$).

\noindent\textit{The above class of problems \eqref{New P} is related to the following well-known problem
 \begin{align}\label{P1}
    \rho \frac{\partial^2u}{\partial t^2} -\left(a+b\int_0^l \bigg|\frac{\partial u}{\partial x}\bigg|^2 \d x \right) \frac{\partial^2 u}{\partial x^2}=0,
\end{align}
which is introduced by Kirchhoff in his celebrated paper \cite{K1883}, as a nonlinear extension of D'Alembert's wave equation for free vibration of elastic strings. Initially \eqref{P1} was studied as the
fundamental problem for understanding several physical, population and biological systems. For further developments in this direction, we refer to \cite{K1883, L1978} and references therein.} One of the earliest studies concerning the existence of a positive, a negative, and a sign-changing solution to the Kirchhoff-type problem involving the Laplacian is due to Zhang and Perera \cite{ZP2006}. In \cite{ZP2006}, the considered the following problem
\begin{equation}\label{P2}
   \begin{aligned}
       \left(a+b\int_\Omega|\nabla u|^2 \d x\right) \left(-\Delta u\right)&=f(u)~\text{ in }\Omega, \\
       u&=0~ \text{ on } \mathbb R^N \setminus \Omega,
       \end{aligned}
   \end{equation}
   where $a,b>0$, $\Omega\subset\mathbb{R}^N$ is a bounded domain with smooth boundary $\partial\Omega$ and the function $f: \mathbb{R}\rightarrow \mathbb{R}$ is locally Lipschitz continuous, subcritical and satisfies some conditions.  Using the invariant sets of descending flows and the Palais-Smale (PS) condition with mountain pass geometry, they obtained the existence of positive, negative, and sign-changing solutions to the problem \eqref{P2}. Later, Mao and Zhang \cite{MZ2009} extended the work of Zhang and Perera \cite{ZP2006}, establishing similar results employing invariant sets of descent flow and minimax methods without using the (PS) condition to $4$-superlinear growth condition on $f$ for $N=1,2,3$. The main difficulty is due to the lack of the Palais–Smale condition. Subsequently, Shuai \cite{S2015} investigated the problem \eqref{P2}, assuming $f\in C^1(\mathbb{R},\mathbb{R})$ satisfying $(f_1)-(f_2)$ and the following Nehari-type monotonicity conditions in place of $(f_3)-(f_4)$ for $N=1,2,3$.
\begin{itemize}
    \item[$\bar{(f_3)}$] $\lim_{|u|\rightarrow \infty}\frac{F(u)}{u^4}=+\infty.$ 
    \item[$(f_4)'$] $\frac{f(u)}{|u|^3}$ is increasing function on $\mathbb{R}\setminus \{0\}$.
\end{itemize}
 \begin{remark}
      Observe that the assumption $(f_4)$ is a weaker than $(f_4)'$ and $\bar{(f_3)}$ is weaker than $(f_3)$. Moreover, the assumptions $(f_4)'$ and $f\in C^1(\mathbb{R},\mathbb{R})$ are significantly stronger than $(f_4)$ and $f\in C(\mathbb{R},\mathbb{R})$. In particular, our nonlinearity is assumed only to be continuous, and no assumptions are made regarding $C^1$. Consequently, our results apply to a wide class of nonlinearities $f$ that are not $C^1$.
      Since many nonlinearities satisfy $(f_4)$ but fail to satisfy $(f_4)'$. For instance, consider the following functions
      $$g_1(t)=|t|^3(2t+|t-1|),~\forall ~t\in \mathbb{R},$$
      and
      % Since many nonlinearities satisfy $(f_4)$ but fail to satisfy $(f_4)'$. 
\begin{align*}\label{Example}
g_2(t)=(\beta|t|^{\frac{3}{2}}+|t|^3)t,~\forall\, t\in \mathbb{R},
\end{align*}
where $0<\beta \leq 4(a\theta \lambda_1)^\frac{1}{2}$ and $\theta\in(0,1)$. Note that $g_1$ is continuous on $\mathbb R$ but not $C^1$ on $\mathbb{R}$ and $g_2$ is satisfies $(f_4)$ but not $(f_4)'$, asserting that our hypothesis admits a wide range of nonlinearities.
  \end{remark}
 In \cite{S2015}, Shuai proved the existence of the least energy sign-changing solution of \eqref{P2} having exactly two nodal domains using the variational tools and deformation lemma. Moreover, the author establishes the convergence of $u_b$ as $b\rightarrow 0$ for a parameter $ b$. For $N=1,2,3$, Tang and Cheng \cite{TC2016} generalized the results of \cite{S2015} when $f$ satisfies $(f_1)-(f_4)$. They established the existence of a least energy sign-changing solution of \eqref{P2} that changes sign only once, using the variational method and introducing a new analytical technique. Moreover, they proved the ``doubling energy property" for sign-changing solutions, using the Non-Nehari manifold method. It is noteworthy to mention that the concept of ``doubling energy" was introduced in Weth \cite{W2006}, which states that the least energy of sign-changing solutions is strictly greater than twice the ``ground-state" energy. Further, Tang and Cheng \cite{TC2016} characterized the nodal domains of least energy sign-changing solutions of \eqref{P2} when $b=0$ and studied the convergent property of $u_b$ as $b\rightarrow 0$. It is almost impossible to enlist all studies of Kirchhoff-type problems involving the Laplacian. The reader may refer to \cite{TC2016, ZP2006, MZ2009, S2015} and the associated references therein for the advancements in this direction.

\par The nonlocal extension of \eqref{P2} involving a general kernel is defined as 
\begin{equation}\label{P3}
   \begin{aligned}
       \left(a+b \int_{\mathbb{R}^{2N}} |u(x)-u(y)|^2K(x-y)\d x\d y\right) \mathcal{L}_Ku&=f(u)~\text{ in }\Omega, \\
       u&=0~ \text{ in } \mathbb R^N \setminus \Omega,
       \end{aligned}
   \end{equation}
where $a,b>0$, $\Omega\subset\mathbb{R}^N$ is a bounded domain with smooth  boundary $\partial\Omega$, $2s<N<4s$, $f$ satisfies $(f_1)-(f_4)$ with the critical exponent $2^*$ being replaced by the fractional Sobolev critical exponent  $2_s^*:=\frac{2N}{N-2s}$ and $\mathcal{L}_K$ is integrodifferential operator with kernel $K$. Observe that $N= 1, 2, 3$ is well-defined under a suitable choice of $s \in (0, 1)$. The existence of least energy sign-changing solutions of \eqref{P3} is guaranteed by Luo \textit{et al.} \cite{LTG2018}. Moreover, they employed a Non-Nehari manifold technique to establish ``doubling energy property" for sign-changing solutions $u_b$. Further, they showed that for any sequence $b_n\searrow 0$ as $n\rightarrow \infty$, there exists a subsequence $b_n$ (still denoted by $b_n$) such that $u_{b_n} \rightarrow u_0$, where $u_0$ is a solution of \eqref{P3} when $b=0$. Subsequently, there have been developments in the direction of the Kirchhoff problem involving the fractional Laplacian. For instance see \cite{G2022, GCF2025, LTG2018} and references therein.

\par Recently, elliptic problems involving mixed local and nonlocal operators have attracted significant interest from researchers for both theoretical aspects and applications to Brownian motion and real-world problems in population dynamics \cite{DV2021} and L\'evy processes \cite{DPV2023}. For the study of regularity and existence of solutions for mixed local-nonlocal operators, we refer to \cite{DV2021, BMV2024, DS2020, DPV2023} and references therein. For instance, in \cite{DV2021, BMV2024}, the authors studied the existence of multiple positive solutions. A first attempt to obtain the regularity of weak solutions has been explored in \cite{GK2022, DM2022}. It is important to emphasize that the majority of existing studies focus on establishing the existence and regularity of positive solutions. In contrast, only a limited number of works address the existence of sign-changing solutions, particularly in the context of problems involving mixed local–nonlocal operators. The first two authors established the existence of sign-changing solutions for the following problem driven by a nonlinear mixed local–nonlocal operator, which is given by
 \begin{equation}\label{P4}
\begin{aligned}
    -\Delta_p u +(-\Delta)^s_p  u&= f(x,u)   \text{ in } \,\,\Omega,\\
      u &=0 \,\text{in} \,~ \mathbb R^N \setminus \Omega,
   \end{aligned}
    \end{equation}
where $\Omega\subset\mathbb{R}^N$ is a bounded domain with smooth boundary $\partial \Omega$, $1<p<\infty$, $0<s<1$, $N>p$, $f:\bar{\Omega}\cross \mathbb{R}\rightarrow \mathbb{R}$ satisfies some conditions and $-\Delta_p u=\div (|\nabla u|^{p-2}\nabla u)$ represents the $p$-Laplace operator and $(-\Delta)^s_p$ denotes the fractional $p$-Laplace operator defined as 
       $$(-\Delta)^s_p u(x) = C({N,s,p}) P.V \int_{\mathbb{R}^N} \frac{(|u(x)-u(y)|^{p-2})(u(x)-u(y))}{|x-y|^{N+ps}} \, \d y.$$
In particular, Bhowmick and Ghosh \cite{BG2025} proved the existence and multiplicity (infinitely many) of sign-changing solutions to \eqref{P4} using invariant sets of descending flows. Moreover, we established the ``doubling energy property" for sign-changing solutions $u_b$, using the Nehari manifold method and degree theory. In \cite{BGS2025}, we established the existence of least energy sign-changing solutions to the problem \eqref{P4} perturbed with a singular term $u^{-\gamma}$. In case of $p=2$, the existence of sign-changing solutions of \eqref{P4} is investigated by Su \textit{et al.} \cite{SVWZ2024}.

% \iffalse
%  Note that for the above study of sign-changing solution of \eqref{P3} when $p=2$, the following conditions are very important.
%  \begin{itemize}
%      \item [$(f_5)$] $\frac{f(u)}{|u|}$ is increasing function on $\mathbb{R}\setminus \{0\}$.
%      \item [$(f_6)$] $\lim_{|u|\rightarrow \infty}\frac{f(u)}{u}=+\infty.$ 
%      \item[$(f_6)'$] There exists $C>0$, such that $F(u)\geq C|u|^\mu$, for some $\mu>2$, where $F(u)=\int_0^u f(s)\d s$.
% \item[$(f_6)''$] There exist $\mu>2$ and $M>0$ such that  $|u|\geq M$,
%    $$f(u)u\geq \mu F(u)>0,~  \text{where}~ F(u):=\int_0^u f(\tau)\d \tau.$$
%  \end{itemize}
%  Observe that, there exists no relation between the conditions of $(f_4)'$ and $(f_5)$, since $g_1(u)=(2u-1)|u|^3$ satisfies $(f_4)'$ but fails to satisfy $(f_5)$, while $g_2(u)=u|u|$ satisfies $(f_5)$ but fails to satisfy $(f_4)'$. \textcolor{red}{Moreover,$(f_6)$ is weak condition of $(f_3)$,} and $(f_6)'$ is stronger condition of $(f_6)$. Finally, in literature, $(f_6)''$ is a well-known (AR) condition, which is a stronger condition than $(f_6)'$. 
% \fi

 \par Motivated by the above-mentioned studies, in this paper, we first consider a Kirchhoff-type problem \eqref{New P} involving mixed local-nonlocal operators with mild assumptions on the nonlinearity $f$. We emphasize that the main difficulty lies in proving that $\mathcal{M}_b \neq 0$, due to the presence of both local and nonlocal operators combined with weak nonlinearities. Moreover, standard approaches like the gradient decomposition method based on the monotonicity of $\frac{f(u)}{|u|^3}$ or the classical $C^1$-manifold method are also not applicable in this setting. To overcome these difficulties, we first establish  Lemma \ref{lmn2.4}, which plays a crucial role in establishing our results. We employ variational methods and the deformation lemma combined with a non-Nehari manifold technique. We mention here that in this paper, we extended the work of Shuai \cite{S2015}, Tang and Cheng \cite{TC2016}, Luo {\it{et al.}} \cite{LTG2018}, and also generalized \cite{BG2025} when $p=2$ and $N=3$ for the Kirchhoff-type problems. To the best of our knowledge, there are no existing results concerning the existence of sign-changing solutions or least energy sign-changing solutions for the Kirchhoff-type problem \eqref{New P}.\\
We now state our main theorems as follows, (for notations and definitions, see Sec. \ref{sec2}).
\begin{theorem}\label{T2.4}
    Assume that $(f_1)-(f_4)$ hold. Then there exists a sign-changing solution $u_b\in \mathcal{M}_b$ of \eqref{New P} such that $\Phi_b(u_b)=m_b$.
\end{theorem}

\begin{theorem}\label{T2.5}
    Assume that $(f_1)-(f_4)$ hold. Then there exists a solution $u\in \mathcal{N}_b$ of \eqref{New P} such that $\Phi_b(u)=c_b>0$. Furthermore, $m_b>2c_b$.
\end{theorem}
\noindent By Theorem \ref{T2.5}, it follows that $m_b>c_b$, which is one of the important theorems in Shuai \cite[Theorem 1.2]{S2015}. Also, observe that the behavior of the problem differs substantially between the cases $b > 0$ and $b = 0$. In the following theorem, we establish a precise relationship between these two cases.
\begin{theorem}\label{T2.6}
    Assume that $(f_1)-(f_4)$ hold. Then there exists a sign-changing solution $v_0\in \mathcal{M}_0$ of \eqref{New P2} such that $\Phi_0(v_0)=m_0$. Moreover, for any sequence $\{b_n\}$ such that $b_n \searrow 0$ as $n \rightarrow \infty$, there exists a subsequence $\{b_n\}$ (still denoted by $\{b_n\}$) such that $u_{b_n}\rightarrow u_0$ in $\mathbb{X}_0^{s,2}(\Omega)$, where $u_0\in \mathcal{M}_0$ is the nodal (sign-changing) solution to the problem \eqref{New P2} such that $\Phi_0(u_0)=m_0$.
\end{theorem}
  \begin{remark}
      It is challenging to provide conclusions regarding the nodal domain of \eqref{New P} and \eqref{New P2} owing to the presence of the nonlocal operator $(-\Delta)^s$. The nonlocal nature of this term significantly complicates the analysis of nodal domains. For further discussions on this issue, refer to Bhowmick, Ghosh, and Saoudi \cite[Remark 1.4]{BGS2025} and  Bhowmick and Ghosh \cite[Remark 1.1]{BG2025}.
  \end{remark}

\par The remaining paper is organized as follows: In Section \ref{sec2}, we recall some fundamental results, important lemmas, and develop the necessary tools for the solution space associated with our problem. In Section \ref{sec3}, we establish the existence of a least energy sign-changing solution to \eqref{New P} using the variational method and deformation lemma. In Section \ref{sec4}, we proved the existence of ground state energy and the ``doubling energy property" for sign-changing solution $u_b$, using a non-Nehari manifold technique. Finally, in Section \ref{sec5}, we prove the existence of least energy solutions to \eqref{New P} for $b=0$ and establish the convergence property of $u_b$ as $b\searrow 0$.

\section{Preliminaries and solution space setup}\label{sec2}

\noindent In this section, we recall some fundamental results of Sobolev spaces and define sign-changing solutions and function spaces associated to it. Let ${ \Omega\subset\mathbb{R}^N}$ be a bounded domain with a smooth boundary $\partial\Omega$. Since we are dealing with the mixed local and nonlocal operators and zero Dirichlet boundary data, it is natural that the solutions are expected to belong to $H^1_0(\Omega)\cap H^s_0(\Omega)$, where $H^1_0(\Omega)$ and $H^s_0(\Omega)$ are standard Sobolev spaces. We now state the following lemma from \cite[Lemma 2.1]{BSM2022} (see also \cite[Proposition 2.2]{NPV2012}), which plays a pivotal role in studying the problem \eqref{New P}.
   \begin{lemma}\label{lm 2.1}
   Let $\Omega\subset\mathbb{R}^N$ be a bounded domain with Lipschitz boundary $\partial\Omega$. Then for $0<s<1$, there exists $C=C(N,s)>0$ such that 
       \begin{equation}\label{eq2.5IM}
       \begin{aligned}
    \|u\|_{H^s(\Omega)} \leq C\|u\|_{H^1(\Omega)}, \,\,\, \forall\, u\in H^1(\Omega).
     \end{aligned}
        \end{equation}
Moreover, for any $u \in H_0^1(\Omega)$ with $u=0$ in $\mathbb R^N\setminus\Omega$, we have
        \begin{equation}\label{sg-emb}
           \int_{\mathbb R^N} \int_{\mathbb R^N} \frac{|u(x)-u(y)|^2}{|x-y|^{N+2s}} \d x \d y \leq C \int_\Omega |\nabla u|^2 \d x 
       \end{equation}
       and in particular, 
        \begin{equation}\label{sg-emb2}
           \int_{\Omega} \int_{\Omega} \frac{|u(x)-u(y)|^2}{|x-y|^{N+2s}} \d x \d y \leq C \int_\Omega |\nabla u|^2 \d x \,\,\, \forall \, u \in H_0^1(\Omega).
           %~\text{with}~u=0~\text{in}~\mathbb R^N\setminus\Omega.
       \end{equation}
       \end{lemma}

\noindent  We now define the solution space corresponding to our problem \eqref{New P}.    
\begin{definition}
  Let $\Omega\subset\mathbb{R}^N$ be a bounded domain with Lipschitz boundary $\partial \Omega$ and let $0<s<1$. We define the space $\mathbb{X}_0^{s,2}(\Omega)$ as the closure of $C_c^\infty(\Omega)$ with respect to the following norm:
  \begin{align}\label{eq2.7 N}
      \|u\|_{\mathbb{X}_0^{s,2}(\Omega)}&=\bigg(\int_\Omega|\nabla u|^2 \d x+\int_{\mathbb{R}^{2N}} \frac{|u(x)-u(y)|^2}{|x-y|^{N+2s}}\d x\d y \bigg)^\frac{1}{2},\,\,\,\,\forall\, u\in C_c^\infty(\Omega).
      \end{align} 
      \end{definition}
    \noindent On using the inequalities \eqref{sg-emb}, \eqref{sg-emb2} and the Poincar\'e inequality in $H_0^1(\Omega)$, we obtain that the norm $\|u\|_{\mathbb{X}_0^{s,2}(\Omega)}=(\int_\Omega|\nabla u|^2 \d x )^\frac{1}{2}$ is equivalent to the norm $$\|u\|_{\mathbb{X}_0^{s,2}(\Omega)}=\bigg(\int_\Omega|\nabla u|^2 \d x+ \int_{\Omega} \int_{\Omega}\frac{|u(x)-u(y)|^2}{|x-y|^{N+2s}}\d x\d y \bigg)^\frac{1}{2}$$ for all $u \in \mathbb{X}_0^{s,2}(\Omega).$ Therefore, the space $\mathbb{X}_0^{s,2}(\Omega)$ is characterized as
    \begin{equation}\nonumber 
          \mathbb{X}_0^{s,2} (\Omega)={\{u\in H_0^1(\Omega): u=0 \text{ in } \mathbb R^N \setminus \Omega} \}.
          \end{equation}
  \noindent Moreover, for $N>2$, using the classical Sobolev inequality, we obtain the following mixed Sobolev inequality in $\mathbb{X}_0^{s,2} (\mathbb{R}^N)$. 
     \begin{equation}\label{eq 2.15}
           \bigg(\int_{\mathbb{R}^N}|u(x)|^{2^*}dx\bigg)^\frac{1}{2^*}\leq C\bigg(\int_{\mathbb{R}^N}|\nabla u(x)|^{2}dx +\int_{\mathbb{R}^{2N}} \frac{|u(x)-u(y)|^2}{|x-y|^{N+2s}}\d x\d y\bigg)^\frac{1}{2},
       \end{equation}
       where $C>0$ is the embedding constant and $2^*:=\frac{2N}{N-2}>\frac{2N}{N-2s}:=2_s^*$. Therefore, we deduce that
          \begin{equation*}
            \|u\|_{L^2{^*}(\Omega)}=\|u\|_{L^2{^*}(\mathbb R^N)} \leq C\|\nabla u\|_{L^2(\mathbb R^N) }\leq C\|u\|_{\mathbb{X}_0^{s,2}(\Omega)},\ \,\,\, \forall\, u \in \mathbb{X}_0^{s,2}(\Omega).
             \end{equation*} 

In particular, the space $\mathbb{X}^{s,2}_0(\Omega)$ is a separable (and reflexive) Hilbert space with respect to the inner product, 
                \begin{equation}
                    \langle u,v\rangle_{\mathbb{X}_0^{s,2}(\Omega)}=\int_{\Omega} \nabla u\cdot\nabla v \d x+\int_{\mathbb{R}^{2N}} \frac{(u(x)-u(y))(v(x)-v(y))}{|x-y|^{N+2s}} \d x \d y,
                \end{equation}
                where ``$\cdot$" denotes the standard scalar product in $\mathbb{R}^N$. Note that the inner product induces the norm \eqref{eq2.7 N} on $\mathbb{X}^{s,2}_0(\Omega)$.
                
                We now state the following embedding results for $\mathbb{X}_0^{s,2}(\Omega)$.
    \begin{theorem}\label{thm cpt}
     Let $s\in (0,1)$ and $\Omega\subset\mathbb{R}^N$ be a bounded domain with Lipschitz boundary $\partial \Omega$. Then, we have the following:
     \begin{itemize}
         \item[(i)] For $N>2$, the embedding $\mathbb{X}_0^{s,2}(\Omega)\hookrightarrow L^q(\Omega)$ is continuous for every $1\leq q\leq 2^*$ and compact for every $1\leq q< 2^*$.
         \item[(ii)] For $N=2$, the embedding $\mathbb{X}_0^{s,2}(\Omega)\hookrightarrow L^q(\Omega)$ is continuous and compact for every $1\leq q< \infty$. Moreover, for $N=1$, the embedding is true even for $q= \infty$.
         %\item For $N<2$, the embedding $\mathbb{X}_0^{s,2}(\Omega)\hookrightarrow L^\infty(\Omega)$ is continuous and compact.
     \end{itemize}
\end{theorem}

Next, we define weak solutions to the problem \eqref{New P}.
\begin{definition}
    We say that $u\in \mathbb{X}_0^{s,2}(\Omega)$ is a weak solution to the problem \eqref{New P} if
    \begin{align}\label{weak form}
        \left(a+b \|u\|_{\mathbb{X}_0^{s,2}(\Omega)}^2\right)  \bigg(\int_\Omega \nabla u\cdot\nabla \phi \d x+&\int_{\mathbb{R}^{2N}} \frac{(u(x)-u(y))(\phi(x)-\phi(y))}{|x-y|^{N+2s}} \d x \d y \bigg) \nonumber \\ 
        = \int_\Omega f(u)\phi \d x, ~\forall~\phi\in \mathbb{X}_0^{s,2}(\Omega).
    \end{align}
\end{definition}
\noindent The associated energy functional $\Phi_b:\mathbb{X}_0^{s,2}(\Omega)\rightarrow \mathbb R$ is defined as
\begin{align}\label{energy}
    {\Phi}_b(u)=\frac{a}{2} \|u\|_{\mathbb{X}_0^{s,2}(\Omega)}^2+\frac{b}{4} \|u\|_{\mathbb{X}_0^{s,2}(\Omega)}^4-\int_\Omega F(u) \d x,
\end{align}
where $F(u)=\int_0^u f(s)\d s$. Clearly, $\Phi_b$ is $C^1$. Therefore, the critical points of $ {\Phi}_b$ are weak solutions to \eqref{New P}. Moreover, we say a weak solution to \eqref{New P} is sign-changing if $u^\pm\neq 0$, where
$u^+=max\{u(x),0\}=\frac{u+|u|}{2}$ and $u^-=min\{u(x),0\}=\frac{u-|u|}{2}.$ 

% \textcolor{red}{We define the set of all nontrivial solutions $\mathcal{N}_b$ and sign-changing solutions $\mathcal{M}_b\subset \mathcal{N}_b$ as}\\
We define the Nehari-type constraint sets $\mathcal{N}_b$ and $\mathcal{M}_b$ as
\begin{align*}
    \mathcal{N}_b&=\{u\in \mathbb{X}_0^{s,2}(\Omega)\setminus \{0\} :\langle\Phi_b'(u),u\rangle=0\}~\text{and}\\
    \mathcal{M}_b&=\{u\in \mathbb{X}_0^{s,2}(\Omega):u^{\pm}\neq0, \langle\Phi_b'(u),u^+\rangle=\langle\Phi_b'(u),u^-\rangle=0\}
\end{align*}
 and set
 $$m_b=\inf_{u\in \mathcal{M}_b} \Phi_b(u)~~\text{and}~~c_b=\inf_{u\in \mathcal{N}_b} \Phi_b(u).$$
 Observe that for any $u=u^++u^-\in {\mathbb{X}_0^{s,2}(\Omega)}$, we get
 \begin{align}
 \|u\|_{\mathbb{X}_0^{s,2}(\Omega)}^2&=\int_\Omega|\nabla u|^2 \d x+\int_{\mathbb{R}^{2N}} \frac{|u(x)-u(y)|^2}{|x-y|^{N+2s}}\d x\d y\nonumber \\
 &= \|u^+\|_{\mathbb{X}_0^{s,2}(\Omega)}^2+ \|u^-\|_{\mathbb{X}_0^{s,2}(\Omega)}^2+2A(u),
 \end{align}
 where \begin{align}
     A(u)=&-\int_{\mathbb{R}^{2N}} \frac{u^-(x)u^+(y)+u^+(x)u^-(y)}{|x-y|^{N+2s}}\d x\d y.  \nonumber
 \end{align}
 Therefore, we have
 \begin{align}
     \Phi_b(u)=&\Phi_b(u^+)+\Phi_b(u^-)+aA(u)+b(A(u))^2+\frac{b}{2}\|u^+\|_{\mathbb{X}_0^{s,2}(\Omega)}^2\|u^-\|_{\mathbb{X}_0^{s,2}(\Omega)}^2\nonumber\\
     &+b\left(\|u^+\|_{\mathbb{X}_0^{s,2}(\Omega)}^2+\|u^-\|_{\mathbb{X}_0^{s,2}(\Omega)}^2\right)A(u), \label{energy 2}\\
     \langle \Phi_b'(u),u\rangle=&\left(a+b \|u\|_{\mathbb{X}_0^{s,2}(\Omega)}^2\right) \|u\|_{\mathbb{X}_0^{s,2}(\Omega)}^2-\int_\Omega f(u)u \d x, \\
     \langle \Phi_b'(u),u^+\rangle=&\langle \Phi_b'(u^+),u^+\rangle+ aA(u)+2b(A(u))^2+b\|u^+\|_{\mathbb{X}_0^{s,2}(\Omega)}^2\|u^-\|_{\mathbb{X}_0^{s,2}(\Omega)}^2 \nonumber \\ &+3b\|u^+\|_{\mathbb{X}_0^{s,2}(\Omega)}^2A(u)+b \|u^-\|_{\mathbb{X}_0^{s,2}(\Omega)}^2 A(u), \\
     \langle \Phi_b'(u),u^-\rangle=&\langle \Phi_b'(u^-),u^-\rangle+ aA(u)+2b(A(u))^2+b\|u^+\|_{\mathbb{X}_0^{s,2}(\Omega)}^2\|u^-\|_{\mathbb{X}_0^{s,2}(\Omega)}^2 \nonumber \\ &+3b\|u^-\|_{\mathbb{X}_0^{s,2}(\Omega)}^2A(u)+b \|u^+\|_{\mathbb{X}_0^{s,2}(\Omega)}^2 A(u).
 \end{align}
We now define weak solutions to the limiting problem \eqref{New P2}.
\begin{definition}
    We say that $u\in \mathbb{X}_0^{s,2}(\Omega)$ is a weak solution to the problem \eqref{New P2} if for all $\phi\in \mathbb{X}_0^{s,2}(\Omega)$, we have
    \begin{align}\label{weak form 2}
       a\bigg(\int_\Omega \nabla u\cdot\nabla \phi \d x+\int_{\mathbb{R}^{2N}} \frac{(u(x)-u(y))(\phi(x)-\phi(y))}{|x-y|^{N+2s}} \d x \d y \bigg)
        = \int_\Omega f(u)\phi \d x.
    \end{align}
\end{definition}
\noindent The associated energy functional $\Phi_0:\mathbb{X}_0^{s,2}(\Omega)\rightarrow \mathbb R$ is $C^1$ which is defined as
\begin{align}\label{energy 0}
    {\Phi}_0(u)=\frac{a}{2} \|u\|_{\mathbb{X}_0^{s,2}(\Omega)}^2-\int_\Omega F(u) \d x,
\end{align}
where $F(u)=\int_0^u f(s)\d s$. Further, we define
\begin{align*}
    \mathcal{N}_0&=\{u\in \mathbb{X}_0^{s,2}(\Omega)\setminus \{0\} :\langle\Phi_0'(u),u\rangle=0\},\\
    \mathcal{M}_0&=\{u\in \mathbb{X}_0^{s,2}(\Omega):u^{\pm}\neq0, \langle\Phi_0'(u),u^+\rangle=\langle\Phi_0'(u),u^-\rangle=0\},\\
    m_0&=\inf_{u\in \mathcal{M}_0} \Phi_0(u)~~\text{and}~~c_0=\inf_{u\in \mathcal{N}_0} \Phi_0(u).
\end{align*}
We conclude this section with the following lemma characterizing the principal eigenvalue and eigenfunction of the mixed local-nonlocal operator \cite{PFR2019} (see also \cite{LGG2024}).
\begin{lemma}
    Let $\Omega$ be an open bounded set of $\mathbb{R}^N$ and $0<s<1$. Consider the following Dirichlet eigenvalue problem:
    \begin{equation}
   \begin{aligned}
        \left(-\Delta u +(-\Delta)^s  u\right)&=\lambda u~\text{ in }\Omega, \\
       u&=0~ \text{ in } \mathbb R^N \setminus \Omega.
       \end{aligned}
   \end{equation}
    Then there exists an eigenvalue $\lambda_1$, which is positive, and it is characterized as follows:
    $$\lambda_1=\min_{u\in \mathbb{X}_0^{s,2}(\Omega)\setminus \{0\}}\frac{\int_\Omega|\nabla u|^2 \d x+\int_{\mathbb{R}^{2N}} \frac{|u(x)-u(y)|^2}{|x-y|^{N+2s}}\d x\d y}{\int_\Omega |u(x)|^2\d x}.$$
\end{lemma}

\section{Existence of least energy sign-changing solutions}\label{sec3}
This section is devoted to guarantee the existence of sign-changing solutions to the problem \eqref{New P}. We begin with the following lemma, which plays a key role in obtaining the main result.
\begin{lemma}\label{lmn2.4}
    Suppose $(f_1)-(f_4)$ hold. For any $k,l\geq 0$ and for all $u\in \mathbb{X}_0^{s,2}(\Omega)$, we have
    \begin{align}
        \Phi_b(u)-\Phi_b(ku^++lu^-)\geq \frac{1-k^4}{4}\langle \Phi_b'(u),u^+\rangle+\frac{1-l^4}{4}\langle \Phi_b'(u),u^-\rangle \nonumber\\
        +\frac{a(1-\theta)(1-k^2)^2}{4}\|u^+\|_{\mathbb{X}_0^{s,2}(\Omega)}^2+\frac{a(1-\theta)(1-l^2)^2}{4}\|u^-\|_{\mathbb{X}_0^{s,2}(\Omega)}^2\nonumber\\
        +\frac{b(k^2-l^2)^2}{4}\|u^+\|_{\mathbb{X}_0^{s,2}(\Omega)}^2\|u^-\|_{\mathbb{X}_0^{s,2}(\Omega)}^2, ~\forall~u\in \mathbb{X}_0^{s,2}(\Omega),
    \end{align}
    where $\theta \in (0,1)$ is given in $(f_4)$.
\end{lemma}
\begin{proof}
    Observe that from $(f_4)$, we get
    \begin{align}\label{f4 inq}
        &\frac{1-l^4}{4}tf(t)+F(lt)-F(t)+\frac{a\theta\lambda_1}{4}(1-l^2)^2 t^2\nonumber \\
        =&\int_l^1\left[\frac{f(t)}{t^3}-\frac{f(kt)}{(kt)^3}+\frac{a\theta\lambda_1(1-k^2)}{(kt)^2} \right] k^3 t^4 \d k\geq 0,~\forall~l\geq 0,~ t\in\mathbb{R}\setminus\{0\}.
    \end{align}
On the other hand, using Young's inequality, $\forall~k,l\geq 0$, we get
\begin{align}\label{young}
    3k^4+l^4-4k^3l \geq& 3k^4+l^4-4\left(\frac{k^4}{\frac{4}{3}}+\frac{l^4}{4}\right)=0 \text{ and } \nonumber\\
    k^4+l^4-4kl+2 \geq& k^4+l^4 -2(k^2+l^2)+2 =(k^2-1)^2+(l^2-1)^2 \geq 0.
\end{align}
    Now from \eqref{energy}, \eqref{f4 inq} and \eqref{young}, we obtain
    \begin{align}
        &\Phi_b(u)-\Phi_b(ku^++lu^-)- \frac{1-k^4}{4}\langle \Phi_b'(u),u^+\rangle-\frac{1-l^4}{4}\langle \Phi_b'(u),u^-\rangle \nonumber \\
        =&\frac{a}{2}\left[\|u^+\|_{\mathbb{X}_0^{s,2}(\Omega)}^2+\|u^-\|_{\mathbb{X}_0^{s,2}(\Omega)}^2+2A(u)-k^2\|u^+\|_{\mathbb{X}_0^{s,2}(\Omega)}^2-l^2\|u^-\|_{\mathbb{X}_0^{s,2}(\Omega)}^2-2klA(u) \right] \nonumber\\
        &+\frac{b}{4}\bigg[\left( \|u^+\|_{\mathbb{X}_0^{s,2}(\Omega)}^2+\|u^-\|_{\mathbb{X}_0^{s,2}(\Omega)}^2+2A(u)\right)^2 \nonumber \\
        &-\left( k^2\|u^+\|_{\mathbb{X}_0^{s,2}(\Omega)}^2+l^2\|u^-\|_{\mathbb{X}_0^{s,2}(\Omega)}^2+2klA(u)\right)^2 \bigg] \nonumber\\
        &+ \left[\int_\Omega F(ku^+)\d x+\int_\Omega F(lu^-)\d x-\int_\Omega F(u^+)\d x-\int_\Omega F(u^-)\d x\right]\nonumber\\
        &- \frac{1-k^4}{4} \bigg[a\left(\|u^+\|_{\mathbb{X}_0^{s,2}(\Omega)}^2+A(u)\right) \nonumber\\
        &+ b\left(\|u^+\|_{\mathbb{X}_0^{s,2}(\Omega)}^2+\|u^-\|_{\mathbb{X}_0^{s,2}(\Omega)}^2+2A(u)\right)\left(\|u^+\|_{\mathbb{X}_0^{s,2}(\Omega)}^2+A(u)\right)-\int_\Omega f(u^+)u^+\d x\bigg]\nonumber\\
        &- \frac{1-l^4}{4} \bigg[a\left(\|u^-\|_{\mathbb{X}_0^{s,2}(\Omega)}^2+A(u)\right) \nonumber\\
        &+ b\left(\|u^+\|_{\mathbb{X}_0^{s,2}(\Omega)}^2+\|u^-\|_{\mathbb{X}_0^{s,2}(\Omega)}^2+2A(u)\right)\left(\|u^-\|_{\mathbb{X}_0^{s,2}(\Omega)}^2+A(u)\right)-\int_\Omega f(u^-)u^-\d x\bigg] \nonumber\\
        =& \frac{a}{4} \|u^+\|_{\mathbb{X}_0^{s,2}(\Omega)}^2+\frac{a}{4}  \|u^-\|_{\mathbb{X}_0^{s,2}(\Omega)}^2-\frac{a}{2}k^2\|u^+\|_{\mathbb{X}_0^{s,2}(\Omega)}^2 -\frac{a}{2}l^2\|u^-\|_{\mathbb{X}_0^{s,2}(\Omega)}^2 +\frac{a}{4}k^4\|u^+\|_{\mathbb{X}_0^{s,2}(\Omega)}^2\nonumber \\
       & +\frac{a}{4}l^4\|u^-\|_{\mathbb{X}_0^{s,2}(\Omega)}^2+\frac{b}{4}(k^2-l^2)^2\|u^-\|_{\mathbb{X}_0^{s,2}(\Omega)}^2\|u^+\|_{\mathbb{X}_0^{s,2}(\Omega)}^2+\frac{b}{2}(k^2-l^2)^2A(u)^2\nonumber\\
       &+ \frac{1}{4}A(u) \bigg[a(2-4kl+k^4+l^4)+b\|u^+\|_{\mathbb{X}_0^{s,2}(\Omega)}^2(3k^4+l^4 -4k^3l) \nonumber \\
       &+b\|u^-\|_{\mathbb{X}_0^{s,2}(\Omega)}^2(3l^4+k^4 -4kl^3) \bigg]+ \int_\Omega \left[\frac{1-k^4}{4} f(u^+)u^+ +F(ku^+)-F(u^+)  \right]\d x \nonumber\\
       &+\int_\Omega \left[\frac{1-l^4}{4} f(u^-)u^-  +F(lu^-)-F(u^-)\right]\d x\nonumber\\
       \geq & \frac{a(1-\theta)(k^2-1)^2}{4}\|u^+\|_{\mathbb{X}_0^{s,2}(\Omega)}^2+\frac{a(1-\theta)(l^2-1)^2}{4}\|u^-\|_{\mathbb{X}_0^{s,2}(\Omega)}^2\nonumber \\
       &+\frac{a\theta(k^2-1)^2}{4} \|u^+\|_{\mathbb{X}_0^{s,2}(\Omega)}^2+\frac{a\theta(l^2-1)^2}{4} \|u^-\|_{\mathbb{X}_0^{s,2}(\Omega)}^2 \nonumber \\
       &+\frac{b}{4}(k^2-l^2)^2\|u^-\|_{\mathbb{X}_0^{s,2}(\Omega)}^2\|u^+\|_{\mathbb{X}_0^{s,2}(\Omega)}^2+ \int_\Omega \left[\frac{1-k^4}{4} f(u^+)u^+ +F(ku^+)-F(u^+)  \right]\d x \nonumber \\
       &+\int_\Omega \left[\frac{1-l^4}{4} f(u^-)u^-  +F(lu^-)-F(u^-)\right]\d x\nonumber\\
       \geq& \frac{a(1-\theta)(k^2-1)^2}{4}\|u^+\|_{\mathbb{X}_0^{s,2}(\Omega)}^2 +\frac{a(1-\theta)(l^2-1)^2}{4}\|u^-\|_{\mathbb{X}_0^{s,2}(\Omega)}^2 \nonumber \\
       &+\frac{b}{4}(k^2-l^2)^2\|u^-\|_{\mathbb{X}_0^{s,2}(\Omega)}^2\|u^+\|_{\mathbb{X}_0^{s,2}(\Omega)}^2\nonumber \\
       &+ \int_\Omega \left[\frac{1-k^4}{4} f(u^+)u^+ +F(ku^+)-F(u^+) +\frac{a\lambda_1\theta(k^2-1)^2}{4} |u^+|^2 \right]\d x \nonumber \\
       &+\int_\Omega \left[\frac{1-l^4}{4} f(u^-)u^-  +F(lu^-)-F(u^-)+\frac{a\lambda_1\theta(l^2-1)^2}{4} |u^-|^2\right]\d x\nonumber \\
       \geq& \frac{a(1-\theta)(k^2-1)^2}{4}\|u^+\|_{\mathbb{X}_0^{s,2}(\Omega)}^2 +\frac{a(1-\theta)(l^2-1)^2}{4}\|u^-\|_{\mathbb{X}_0^{s,2}(\Omega)}^2 \nonumber \\
       &+\frac{b}{4}(k^2-l^2)^2\|u^-\|_{\mathbb{X}_0^{s,2}(\Omega)}^2\|u^+\|_{\mathbb{X}_0^{s,2}(\Omega)}^2.\nonumber
    \end{align}
    This completes the proof.
\end{proof}
\begin{cor}\label{cor}
   Let $(f_1)-(f_4)$ hold. For any $k,l\geq 0$, we have
    \begin{align}
        &\Phi_b(u)-\Phi_b(ku^++lu^-)\geq 
        \frac{a(1-\theta)(1-k^2)^2}{4}\|u^+\|_{\mathbb{X}_0^{s,2}(\Omega)}^2  \nonumber\\
        &+\frac{a(1-\theta)(1-l^2)^2}{4}\|u^-\|_{\mathbb{X}_0^{s,2}(\Omega)}^2
        +\frac{b(k^2-l^2)^2}{4}\|u^+\|_{\mathbb{X}_0^{s,2}(\Omega)}^2\|u^-\|_{\mathbb{X}_0^{s,2}(\Omega)}^2, ~\forall~u\in \mathcal{M}_b,
    \end{align}
    where $\theta \in (0,1)$ is given $(f_4)$. Moreover,
    $$ \Phi_b(u)=\max_{k,l \geq 0}\Phi_b(ku^++lu^-) ~\forall~u\in \mathcal{M}_b.$$
\end{cor}

 Next, we prove that $\mathcal{M}_b\neq 0$ for every $u\in \mathbb{X}_0^{s,2}(\Omega)$ with $u^{\pm}\neq0$.

\begin{lemma}\label{lmn5.4}
    Let the assumptions $(f_1)-(f_4)$ hold. If $u\in \mathbb{X}_0^{s,2}(\Omega)$ with $u^{\pm}\neq0$, then there exists a unique pair $(k_u,l_u)$ of positive numbers such that 
    $k_uu^++l_uu^-\in\mathcal{M}_b.$
\end{lemma}
\begin{proof}
We divide the proof into two parts.\\

 \noindent\underline{\textbf{{Existence part:}}} For $k,l>0$, let us consider two functions, $h_1$ and $h_2$, which are defined as follows:
    \begin{align}\label{eq5.150}
    & h_1(k,l)=\langle \Phi'(ku^++lu^-), ku^+\rangle  \nonumber \\
    =&\left(a+b \|ku^++lu^-\|_{\mathbb{X}_0^{s,2}(\Omega)}^2\right)  \bigg(\int_\Omega \nabla (ku^++lu^-)\cdot\nabla (ku^+) \d x\nonumber \\
    &+\int_{\mathbb{R}^{2N}} \frac{\{(ku^++lu^-)(x)-(ku^++lu^-)(y)\}\{(ku^+)(x)-(ku^+)(y)\}}{|x-y|^{N+2s}} \d x \d y \bigg) \nonumber \\
    &-\int_\Omega f(ku^++lu^-)ku^+ \d x \nonumber \\
    =& \left(a+b k^2 \|u^+\|_{\mathbb{X}_0^{s,2}(\Omega)+}^2+b l^2 \|u^-\|_{\mathbb{X}_0^{s,2}(\Omega)}^2+2bklA(x)\right) \left(k^2 \|u^+\|_{\mathbb{X}_0^{s,2}(\Omega)+}^2 +lk A(u)\right)  \nonumber \\
    &-\int_\Omega f(ku^+)ku^+ \d x \nonumber \\
    =& ak^2 \|u^+\|_{\mathbb{X}_0^{s,2}(\Omega)+}^2+akl A(u)+ bk^4\|u^+\|_{\mathbb{X}_0^{s,2}(\Omega)}^4+bk^2l^2\|u^+\|_{\mathbb{X}_0^{s,2}(\Omega)}^2\|u^+\|_{\mathbb{X}_0^{s,2}(\Omega)}^2 \nonumber \\
    &+3bk^3l \|u^+\|_{\mathbb{X}_0^{s,2}(\Omega)}^2 A(u)+bkl^3 \|u^-\|_{\mathbb{X}_0^{s,2}(\Omega)}^2A(u)+2bk^2l^2A(u)^2\nonumber\\
    &-\int_\Omega f(ku^+)ku^+ \d x,
    \end{align}
    and 
    \begin{align}\label{eq5.160}
       & h_2(k,l)=\langle \Phi'(ku^++lu^-), lu^-\rangle \nonumber \\
       =& al^2 \|u^-\|_{\mathbb{X}_0^{s,2}(\Omega)+}^2+akl A(u)+ bl^4\|u^-\|_{\mathbb{X}_0^{s,2}(\Omega)}^4+bk^2l^2\|u^+\|_{\mathbb{X}_0^{s,2}(\Omega)}^2\|u^+\|_{\mathbb{X}_0^{s,2}(\Omega)}^2 \nonumber \\
    &+3bkl^3 \|u^-\|_{\mathbb{X}_0^{s,2}(\Omega)}^2 A(u)+bk^3l \|u^+\|_{\mathbb{X}_0^{s,2}(\Omega)}^2A(u)+2bk^2l^2A(u)^2\nonumber\\
    &-\int_\Omega f(lu^-)lu^- \d x.
    \end{align}
   From $(f_1)-(f_3)$, we deduce that $h_1(k,k)>0$ and $h_2(k,k)>0$ for small $k>0$, whereas $h_1(l,l)<0$ and $h_2(l,l)<0$ for large $l>0$. Therefore, there exist $r>0$ small enough and $R>0$ sufficiently large such that 
\begin{align}\label{eq5.180}
   & h_1(r,r)>0,~h_2(r,r)>0,\\
    \text{and}~& h_1(R,R)<0,~h_2(R,R)<0, \text{ with }0<r<R.\label{eq5.190}
\end{align}
On using \eqref{eq5.150}, \eqref{eq5.160}, \eqref{eq5.180} and \eqref{eq5.190}, we get
\begin{align}\label{eq5.200}
   & h_1(r,l)\geq h_1(r,r)>0,~h_1(R,l)\leq h_1(R,R)<0,~\forall~l\in[r,R],\\
    \text{and}~& h_2(k,r)\geq h_2(r,r)>0,~h_2(k,R)\leq h_2(R,R)<0,~\forall~k\in(r,R].\label{eq5.210}
\end{align}
    Thus, applying Miranda's theorem \cite{M1940}, we obtain that there exist $k_u,l_u\in [r,R]$ such that $h_1(k_u,l_u)=0$ and $h_2(k_u,l_u)=0.$ Therefore, we get $k_uu^++l_uu^-\in \mathcal{M}_b$.\\
    
\noindent \underline{\textbf {Uniqueness part:}} Let $(k_1,l_1)$ and $(k_2,l_2)$ represent two distinct positive pairs such that $k_iu^++l_iu^-\in \mathcal{M}_b$ for $i=1,2$. Using Corollary \ref{cor}, we obtain
\begin{align}\label{15}
        \Phi_b(k_1u^++l_1u^-)-\Phi_b(k_2u^++l_2u^-)\geq 
        +\frac{a(1-\theta)(k_1^2-k_2^2)^2}{4{k_1^2}}\|u^+\|_{\mathbb{X}_0^{s,2}(\Omega)}^2  \nonumber\\
        +\frac{a(1-\theta)(l_1^2-l_2^2)^2}{4{l_1^2}}\|u^-\|_{\mathbb{X}_0^{s,2}(\Omega)}^2,
    \end{align}
    and
    \begin{align}\label{16}
        \Phi_b(k_2u^++l_2u^-)-\Phi_b(k_1u^++l_1u^-)\geq 
        +\frac{a(1-\theta)(k_2-k_1^2)^2}{4{k_2^2}}\|u^+\|_{\mathbb{X}_0^{s,2}(\Omega)}^2  \nonumber\\
        +\frac{a(1-\theta)(l_2^2-l_1^2)^2}{4{l_2^2}}\|u^-\|_{\mathbb{X}_0^{s,2}(\Omega)}^2,
    \end{align}
    Therefore, from \eqref{15} and \eqref{16}, we conclude that $k_1=k_2$ and $l_1=l_2$.
\end{proof}

% $\Phi_b(u)=\max_{s,t \geq 0}\Phi_b(ku^++lu^-)$
\noindent As a consequence of Corollary \ref{cor} and Lemma \ref{lmn5.4}, we have following lemma.
\begin{lemma}\label{lmn5.5}
    Let the assumptions $(f_1)-(f_4)$ hold. Then, we have
    \begin{equation}\label{gg3.13}
        \inf_{u\in \mathcal{M}_b}\Phi_b(u)=m_b=\inf_{u\in {\mathbb{X}_0^{s,2}(\Omega)}, u^\pm \neq 0} \max_{k,l \geq 0}\Phi_b(ku^++lu^-).
    \end{equation} 
\end{lemma}
\begin{proof}
    From Corollary \ref{cor}, we obtain
    \begin{align}\label{eq17}
        \inf_{u\in {\mathbb{X}_0^{s,2}(\Omega)}, u^\pm \neq 0} \max_{k,l \geq 0}\Phi_b(ku^++lu^-) &\leq \inf_{u\in \mathcal{M}_b}\max_{k,l \geq 0}\Phi_b(ku^++lu^-)\nonumber\\
        &=\inf_{u\in \mathcal{M}_b}\Phi_b(u)=m_b.
    \end{align}
    On the other side, using Lemma \ref{lmn5.4}, there exists a unique pair $(k_u,l_u)$ of positive numbers such that 
    $k_uu^++l_uu^-\in\mathcal{M}_b$ for any $u\in \mathbb{X}_0^{s,2}(\Omega)$ with $u^{\pm}\neq0$. Thus, we get
    \begin{align}\label{eq18}
        \max_{k,l \geq 0}\Phi_b(ku^++lu^-) \geq \Phi_b(k_uu^++l_uu^-)\geq \inf_{w\in \mathcal{M}_b}\Phi_b(w)=m_b.
    \end{align}
    From \eqref{eq18}, we have
    \begin{align}\label{eq19}
        \inf_{u\in {\mathbb{X}_0^{s,2}(\Omega)}, u^\pm \neq 0} \max_{k,l \geq 0}\Phi_b(ku^++lu^-) \geq  \inf_{u\in \mathcal{M}_b}\Phi_b(u)=m_b.
    \end{align}
    Therefore, combining \eqref{eq17} and \eqref{eq19}, we obtain \eqref{gg3.13}, competing the proof.
\end{proof}
The following lemma guarantees the existence of $u_b\in\mathcal{M}_b$ such that $\Phi_b(u_b)=m_b$.
\begin{lemma}\label{lmn5.6}
     Let the assumptions $(f_1)-(f_4)$ hold. Then $m_b=\inf_{u\in \mathcal{M}_b}\Phi_b(u)>0$ is achieved for some $u_b\in\mathcal{M}_b$.
\end{lemma}
\begin{proof}
    Recall that, for $l=0$ in \eqref{f4 inq}, we get
        \begin{align}\label{eq21}
            \frac{1}{4}tf(t)-F(t)+\frac{a\theta\lambda_1}{4} t^2\geq 0,~ t\in\mathbb{R}.
        \end{align}
        From conditions $(f_1)-(f_2)$, we have for any $\epsilon>0$, there exists $C_\epsilon>0$ such that
    \begin{align}\label{eq20}
        f(x,t)t\leq \epsilon |t|^2+C_\epsilon |t|^q,~\forall~t\in \mathbb{R}.
        \end{align}
        \underline{\textbf{Claim $m_b>0$ and $m_b$ is achieved:}}\\
        
       \noindent \textit{Proof of the Claim:} Since, $u\in \mathcal{M}_b$, we have $\langle \Phi_b(u),u \rangle=0$. Now, using \eqref{eq20} and Sobolev embedding Theorem \ref{thm cpt}, we obtain
        \begin{align}\label{eq22}
            a\|u\|_{\mathbb{X}_0^{s,2}(\Omega)}^2\leq a\|u\|_{\mathbb{X}_0^{s,2}(\Omega)}^2+b\|u\|_{{\mathbb{X}_0^{s,2}(\Omega)}}^4=&\int_\Omega f(u) u \d x  \nonumber \\
            \leq & \int_\Omega \left(\epsilon |u|^2+C_\epsilon |u|^q \right) \d x\nonumber\\
             \leq &  \epsilon C \|u\|^2_{\mathbb{X}_0^{s,2}(\Omega)}+C_pC_\epsilon \|u\|^q_{\mathbb{X}_0^{s,2}(\Omega)}
        \end{align}
        for all $u\in  \mathcal{M}_b$, where $C$ and $C_p$ are embedding constants. Choosing $\epsilon=\frac{a}{2C}$, we get 
        $$\|u\|_{\mathbb{X}_0^{s,2}(\Omega)}^2 \geq \left(\frac{a}{2C_pC_\epsilon} \right)^\frac{2}{(q-2)}:=k >0.$$
        Using \eqref{energy}, \eqref{eq21} and \eqref{eq20}, we obtain
        \begin{align}\label{eq23}
            \Phi_b(u)=&\Phi_b(u) -\frac{1}{4}\langle \Phi'_b(u),u \rangle=\frac{a}{4}\|u\|_{\mathbb{X}_0^{s,2}(\Omega)}^2+\int_\Omega \left(\frac{1}{4}f(u)u-F(u) \right) \d x \nonumber \\
            \geq &\frac{a}{4}\|u\|_{\mathbb{X}_0^{s,2}(\Omega)}^2-\frac{a\theta\lambda_1}{4} \int_\Omega u^2\d x \geq \frac{a(1-\theta)}{4} \|u\|_{\mathbb{X}_0^{s,2}(\Omega)}^2 \geq \frac{ak(1-\theta)}{4}.
        \end{align}
        Thus, from \eqref{eq23}, we conclude that $m_b\geq \frac{ak(1-\theta)}{4} >0$, completing the first part.\\
        
        \par For the second part, let us $\{u_n\}\subset \mathcal{M}_b$ be a sequence such that $\Phi_b(u_n)\rightarrow m_b$ as $n\rightarrow \infty$. Then for large $n$, we get
        \begin{align}\label{eq24}
            m_b+1 \geq \Phi_b(u_n) -\frac{1}{4}\langle \Phi'_b(u_n),u_n \rangle \geq \frac{a(1-\theta)}{4} \|u_n\|_{\mathbb{X}_0^{s,2}(\Omega)}^2.
        \end{align}
        Thus, from \eqref{eq24}, we conclude that $\{u_n\}$ is bounded in $\mathbb{X}_0^{s,2}(\Omega)$. Therefore, there exists a subsequence $\{u_n\}$ (still denoted by $\{u_n\}$) and $u_b\in \mathbb{X}_0^{s,2}(\Omega)$ such that 
        \begin{align}
            u_n^\pm &\rightharpoonup u_b^\pm\in \mathbb{X}_0^{s,2}(\Omega)~\text{weakly in}~ \mathbb{X}_0^{s,2}(\Omega),\\
        u_n^\pm &\rightarrow u_b^\pm~\text{strongly in}~ L^r(\Omega)~\text{for}~r\in[1,2^*), \label{strong}\\ 
        u_n(x)& \rightarrow u_b(x)~a.e.~\text{in}~\Omega.
        \end{align}
        Moreover, since $u_n \in \mathcal{M}_b$, we have $\langle \Phi'_b(u_n),u_n^\pm \rangle=0$, that is,
        \begin{align}\label{25}
            \left(a+b\|u_n\|^2_{\mathbb{X}_0^{s,2}(\Omega)}\right)\left(\|u_n^\pm\|^2_{\mathbb{X}_0^{s,2}(\Omega)}+A(u_n) \right)=\int_\Omega f(u_n^\pm)u^\pm_n \d x. 
            \end{align}
            We split the remaining part of the proof into three parts. We first prove $ u_b^\pm \neq 0$, in the next part, we show $\langle \Phi'_b(u_b),u_b^\pm \rangle \leq0$. Finally, we prove that $\Phi_b(u_b)=m_b$ and $u_b\in \mathcal{M}_b$.\\

           \noindent On using \eqref{25} with a similar argument as in \eqref{eq22}, we get
        \begin{align}\label{eq3.26}
            a\|u_n^\pm\|^2_{\mathbb{X}_0^{s,2}(\Omega)} \leq \int_\Omega f(u_n^\pm)u^\pm_n \d x\leq  \epsilon C \|u^\pm\|^2_{\mathbb{X}_0^{s,2}(\Omega)}+C_pC_\epsilon \|u^\pm\|^q_{\mathbb{X}_0^{s,2}(\Omega)}.
        \end{align}
        Thus, \eqref{eq3.26} with a suitable $\epsilon$, there exists $k_1>0$ such that $\|u_n^\pm\|^2_{\mathbb{X}_0^{s,2}(\Omega)} \geq k_1$. On the other hand, we have
        \begin{align}
            ak_1 \leq a\|u_n^\pm\|^2_{\mathbb{X}_0^{s,2}(\Omega)} \leq \int_\Omega f(u_n^\pm)u^\pm_n \d x \leq \int_\Omega \left(\epsilon |u_n^\pm|^2+C_\epsilon |u_n^\pm|^q \right) \d x.
        \end{align}
        Using the boundedness of $\{u_n\}$ in $\mathbb{X}_0^{s,2}(\Omega)$ and Sobolev embedding Theorem \ref{thm cpt}, there exists a $C>0$ such that
        \begin{align}
            ak_1 \leq \epsilon C +C_\epsilon \int_\Omega |u_n^\pm|^q \d x.
        \end{align}
        Now, choosing $\epsilon=\frac{ak_1}{2C}$, we obtain
        \begin{align}
            \int_\Omega |u_n^\pm|^q \d x \geq \frac{ak_1}{2C_\epsilon}:=k_2>0.
        \end{align}
         Therefore, using \eqref{strong}, we get
          \begin{align}
            \int_\Omega |u_b^\pm|^q \d x \geq \frac{ak_1}{2C_\epsilon}=k_2.
        \end{align}
       Thus, we have $u_b^\pm\neq0$.\\

        \noindent Following the arguments as in \cite{W1996}, combined with the assumptions $(f_1)-(f_3)$, we deduce that
    \begin{align} \label{eq25}
        &\lim_{n\rightarrow \infty} \int_\Omega f(u_n^\pm)u_n^\pm \d x =\int_\Omega f(u_b^\pm)u_b^\pm \d x, ~\text{and}\\ \label{eq26}
        &\lim_{n\rightarrow \infty} \int_\Omega F(u_n^\pm) \d x =\int_\Omega F(u_b^\pm) \d x.
    \end{align}
    Since $A(u_n)\geq 0$, we use Fatou's lemma to obtain
    \begin{align}\label{eq27}
        A(u_b) \leq \liminf_{n\rightarrow \infty} A(u_n).
    \end{align}
    Applying the weak semicontinuity of norm and \eqref{eq25}-\eqref{eq27}, we get
    \begin{align}
        &\left(a+b\|u_b\|^2_{\mathbb{X}_0^{s,2}(\Omega)}\right)\left(\|u_b^\pm\|^2_{\mathbb{X}_0^{s,2}(\Omega)}+A(u_b) \right) \nonumber \\
        \leq &\liminf_{n\rightarrow \infty} \left\{\left(a+b\|u_n\|^2_{\mathbb{X}_0^{s,2}(\Omega)}\right)\left(\|u_n^\pm\|^2_{\mathbb{X}_0^{s,2}(\Omega)}+A(u_n) \right) \right\} \nonumber \\
        = &\liminf_{n\rightarrow \infty} \int_\Omega f(u_n^\pm)u^\pm_n \d x\leq \int_\Omega f(u_b^\pm)u^\pm_b \d x.
    \end{align}
    Therefore, we conclude that $\langle \Phi'_b(u_b),u_b^\pm \rangle \leq0$.\\

    \noindent Finally, using \eqref{energy}, \eqref{eq21}, Lemma \ref{lmn2.4}, the weak semicontinuity of norm and Fatou's lemma, we get
    \begin{align}
        m_b=&\lim_{n\rightarrow\infty} \bigg[\Phi_b(u_n)-\frac{1}{4}\langle\Phi_b'(u_n),u_n\rangle\bigg]\nonumber \\
        =&\lim_{n\rightarrow\infty} \bigg[\frac{a}{4}\|u_n\|_{\mathbb{X}_0^{s,2}(\Omega)}^2+\int_\Omega \left(\frac{1}{4}f(u_n)u_n-F(u_n) \right) \d x \bigg]\nonumber \\
        \geq & \frac{a}{4} \liminf_{n\rightarrow\infty} \left(\|u_n\|_{\mathbb{X}_0^{s,2}(\Omega)}^2-\theta\lambda_1\|u_n\|^2_{L^2(\Omega)} \right)\nonumber\\
        &+\liminf_{n\rightarrow\infty}\int_\Omega \left(\frac{1}{4}f(u_n)u_n-F(u_n) +\frac{a\theta\lambda_1}{4}|u_n|^2\right) \d x \nonumber \\
        \geq &\frac{a}{4}  \left(\|u_b\|_{\mathbb{X}_0^{s,2}(\Omega)}^2-\theta\lambda_1\|u_b\|^2_{L^2(\Omega)} \right)+\int_\Omega \left(\frac{1}{4}f(u_b)u_b-F(u_b) +\frac{a\theta\lambda_1}{4}|u_b|^2\right) \d x \nonumber \\
        = &\frac{a}{4}  \|u_b\|_{\mathbb{X}_0^{s,2}(\Omega)}^2+\int_\Omega \left(\frac{1}{4}f(u_b)u_b-F(u_b) \right) \d x=\Phi_b(u_b)-\frac{1}{4}\langle\Phi_b'(u_b),u_b\rangle \nonumber \\
        \geq & \sup_{k,l\geq 0}\left( \Phi_b(ku_b^++lu_b^-) +\frac{1-k^4}{4}\langle \Phi_b'(u_b),u_b^+\rangle+\frac{1-l^4}{4}\langle \Phi_b'(u_b),u_b^-\rangle \right)  \nonumber \\
        &-\frac{1}{4}\langle\Phi'_b(u_b),u_b\rangle \nonumber \\
        \geq & \sup_{k,l\geq 0}\left( \Phi_b(ku_b^++lu_b^-)-\frac{k^4}{4}\langle \Phi_b'(u_b),u_b^+\rangle-\frac{l^4}{4}\langle \Phi_b'(u_b),u_b^-\rangle \right)\nonumber \\
         \geq & \sup_{k,l\geq 0}\Phi_b(ku_b^++lu_b^-)\geq m_b.
    \end{align}
    Thus, we get
    $$\lim_{n\rightarrow \infty}\left(\|u_n\|_{\mathbb{X}_0^{s,2}(\Omega)}^2-\theta\lambda_1\|u_n\|^2_{L^2(\Omega)} \right)=\left(\|u_b\|_{\mathbb{X}_0^{s,2}(\Omega)}^2-\theta\lambda_1\|u_b\|^2_{L^2(\Omega)} \right).$$
    Therefore, $u_n\rightarrow u_b$ in ${\mathbb{X}_0^{s,2}(\Omega)}$, implying that $\Phi_b(u_b)=m_b$ and $u_b\in \mathcal{M}_b$. This completes the proof.
\end{proof}
In the next lemma, we use the deformation Lemma to prove that $\Phi_b'(u_1)=0$ when $\Phi_b(u_1)=m_b$ for some $u_1\in \mathcal{M}_b$.
\begin{lemma}\label{lmn5.7}
    If $\Phi_b(u_1)=m_b$ for some $u_1\in \mathcal{M}_b$, then $u_1$ is a critical point of $\Phi_b$.
\end{lemma}
\begin{proof}
    We proceed with the method of contradiction. Let $\Phi'_b(u_1)\neq0$, then there exist $\rho,\mu>0$ such that
    $$\|\Phi'_b(u)\|_{\mathbb{X}_0^{s,2}(\Omega)}\geq\rho,~\forall ~u\in B_{3\mu}(u_1),$$
    where $B_{3\mu}(u_1)=\{u\in{\mathbb{X}_0^{s,2}(\Omega)}:\|u-u_1\|_{\mathbb{X}_0^{s,2}(\Omega)}\leq3\mu\}$. Moreover, we have
    \begin{align}\label{eq1}
        \Phi_b(ku_1^++lu_1^-)\leq &\Phi_b(u_1)-\frac{a(1-\theta)(1-k^2)^2}{4}\|u_1^+\|_{\mathbb{X}_0^{s,2}(\Omega)}^2 
        -\frac{a(1-\theta)(1-l^2)^2}{4}\|u_1^-\|_{\mathbb{X}_0^{s,2}(\Omega)}^2\nonumber \\
        =&m_b-\frac{a(1-\theta)(1-k^2)^2}{4}\|u_1^+\|_{\mathbb{X}_0^{s,2}(\Omega)}^2 
        -\frac{a(1-\theta)(1-l^2)^2}{4}\|u_1^-\|_{\mathbb{X}_0^{s,2}(\Omega)}^2,
    \end{align}
    for all $k,l\geq 0$. Let $R:=\left(\frac{1}{2},\frac{3}{2}\right)\cross\left(\frac{1}{2},\frac{3}{2}\right)$. Using \eqref{eq1}, we obtain
    \begin{align}\label{eq2}
        \delta:= \max_{(k,l)\in \partial R} \Phi_b(ku_1^++lu_1^-) < m_b.
    \end{align}
    Choose $\epsilon=\min\{\frac{m_b-\delta}{3},1,\frac{\rho\mu}{8}\}$ and $B_\mu(u_1)=\{u\in{\mathbb{X}_0^{s,2}(\Omega)}:||u-u_1||_{\mathbb{X}_0^{s,2}(\Omega)}\leq\mu_1\}$. Then, using the quantitative deformation Lemma 2.3 in \cite{W1996}, there exists a continuous mapping $\eta:[0,1]\cross {\mathbb{X}_0^{s,2}(\Omega)} \rightarrow {\mathbb{X}_0^{s,2}(\Omega)}$ such that \\
    $(a)$ $\eta(1,u)=u$ if $u\notin \Phi_b^{-1}[m_b-2\epsilon,m_b+2\epsilon_1]\cap B_{2\mu}(u_1)$,\\
    $(b)$ $\eta(1,\Phi^{m_b+\epsilon}\cap B_{\mu}(u_1))\subset \Phi_b^{m_b-\epsilon}$,\\
    $(c)$ $\Phi_b(\eta(1,u))\leq \Phi_b(u)$, $\forall~u \in {\mathbb{X}_0^{s,2}(\Omega)}$.\\
    Now, applying Corollary \eqref{cor}, we obtain $\Phi_b(ku_1^++lu_1^-) \leq \Phi_b(u_1)=m_b$ for all $k,l \geq0$. Then using $(b)$, we get
    \begin{align}\label{eq4}
        \Phi_b(\eta(1,ku_1^++lu_1^-)) \leq m_b-\epsilon ~\forall~ k,l \geq 0 \text{ and } |k-1|^2+|k-1|^2 < \frac{\mu^2}{2\|u_1\|^2_{\mathbb{X}_0^{s,2}(\Omega)}}.
    \end{align}    
    On the other hand, using $(c)$ and \eqref{eq1}, we get
    \begin{align}
        \Phi_b(\eta(1,ku_1^++lu_1^-)) \leq &\Phi_b(ku_1^++lu_1^-) \nonumber \\
        \leq & m_b-\frac{a(1-\theta)(1-k^2)^2}{4}\|u_1^+\|_{\mathbb{X}_0^{s,2}(\Omega)}^2 
        -\frac{a(1-\theta)(1-l^2)^2}{4}\|u_1^-\|_{\mathbb{X}_0^{s,2}(\Omega)}^2 \nonumber \\
        \leq & m_b- \frac{a(1-\theta)\mu^2}{8\|u_1\|^2_{\mathbb{X}_0^{s,2}(\Omega)}} \min\{\|u_1^+\|^2_{\mathbb{X}_0^{s,2}(\Omega)}, \|u_1^-\|^2_{\mathbb{X}_0^{s,2}(\Omega)}\}, \label{eq5}
    \end{align}
    for all $ k,l \geq 0$ and $|k-1|^2+|l-1|^2 \geq \frac{\mu^2}{2\|u_1\|^2_{\mathbb{X}_0^{s,2}(\Omega)}}$. Therefore, from \eqref{eq4} and \eqref{eq5}, we deduce that
    \begin{align}\label{eq6}
        \max_{(k,l)\in \bar{R}}\Phi_b(\eta(1,ku_1^++lu_1^-)) <m_b.
    \end{align}
    The remaining part of the proof follows from \cite{S2015} (see also \cite[Lemma 2.13]{TC2017}).
\end{proof}

\noindent We now prove our first main theorem.
\section*{{\bf\emph{Proof of the Theorem \ref{T2.4}}}}
%We divide the proof into two parts.\\
%\textbf{Part-I: (Existence).} 
\begin{proof}
From Lemma \ref{lmn5.6} and Lemma \ref{lmn5.7}, there exists $u_b\in \mathcal{M}_b$ such that $\Phi_b(u_b)=m_b$ and $\Phi'_b(u_b)=0$. Therefore, $u_b$ is a sign-changing solution to the problem \eqref{New P}.
\end{proof}

\section{Existence of ground state solutions}\label{sec4}
In this section, we use the non-Nehari manifold approach (Introduced by Tang \cite{T2014}) to prove the existence of a ground state solution to the problem \eqref{New P}. We proceed with proving the following important lemma.
\begin{lemma}\label{lmn1}
     Let the assumptions $(f_1)-(f_4)$ hold. Then for any $k\geq 0$, we have
    \begin{align}
        \Phi_b(u)-\Phi_b(ku)\geq \frac{1-k^4}{4}\langle \Phi_b'(u),u\rangle
        +\frac{a(1-\theta)(1-k^2)^2}{4} \|u\|_{\mathbb{X}_0^{s,2}(\Omega)}^2 ~\forall~u\in \mathbb{X}_0^{s,2}(\Omega),
    \end{align}
    where $\theta \in (0,1)$ is given in $(f_4)$.
\end{lemma}
\begin{proof}
    The inequality follows immediately from \eqref{energy} and \eqref{f4 inq}. Indeed,
    \begin{align}
       & \Phi_b(u)-\Phi_b(ku)- \frac{1-k^4}{4}\langle \Phi_b'(u),u\rangle
        -\frac{a(1-\theta)(1-k^2)^2}{4} \|u\|_{\mathbb{X}_0^{s,2}(\Omega)}^2 \nonumber \\
        = & \int_\Omega \left(\frac{1-k^4}{4} f(u)u +F(ku)-F(u) \right) \d x+\frac{a \theta(1-k^2)^2}{4} \|u\|_{\mathbb{X}_0^{s,2}(\Omega)}^2  \nonumber \\
        \geq & \int_\Omega \left(\frac{1-k^4}{4} f(u)u +F(ku)-F(u)+\frac{a \lambda_1\theta(1-k^2)^2}{4} |u|^2\right) \d x \geq 0.\nonumber
    \end{align}
    This completes the proof.
\end{proof}
\begin{cor} \label{cor2}
    Suppose $(f_1)-(f_4)$ hold. For any $k\geq 0$, we have
    \begin{align}
        \Phi_b(u)\geq \Phi_b(ku)
        +\frac{a(1-\theta)(1-k^2)^2}{4}\|u\|_{\mathbb{X}_0^{s,2}(\Omega)}^2 ~\forall~u\in \mathcal{N}_b,
    \end{align}
    where $\theta \in (0,1)$ is given in $(f_4)$. Moreover,
    $$ \Phi_b(u)=\max_{k\geq 0}\Phi_b(ku).$$
\end{cor}
\begin{lemma}\label{lmn2}
     Let the assumptions $(f_1)-(f_4)$ hold. If $u\in \mathbb{X}_0^{s,2}(\Omega)\setminus \{0\}$, then there exists a unique $k_u>0$ such that 
    $k_uu\in\mathcal{N}_b.$
\end{lemma}
\begin{proof}
    The proof follows from a similar argument as in the proof of Lemma \ref{lmn5.4} with the function $h(k)=\langle \Phi'(ku), ku\rangle$ and Corollary \ref{cor2}.
\end{proof}

\begin{lemma}\label{lmn3}
    Suppose that $(f_1)-(f_4)$ hold. Then 
    $$\inf_{u\in \mathcal{N}_b}\Phi_b(u)=c_b=\inf_{u\in {\mathbb{X}_0^{s,2}(\Omega)}, u\neq 0} \max_{k,l \geq 0}\Phi_b(ku) >0.$$ 
\end{lemma}
\begin{proof}
The proof is immediate by following similar arguments as in the proof of Lemma \ref{lmn5.5} combined with Corollary \ref{cor2} and Lemma \ref{lmn2}. The proof of $c_b>0$ follows similarly from the proof of $m_b>0$ in Lemma \ref{lmn5.6}.
\end{proof}
%In the following lemma, we prove the existence of the sequence $u_n \in \mathbb{X}_0^{s,2}(\Omega)$, which is not directly chosen from the Nehari manifold of $\Phi_b$.
\begin{lemma}\label{lmn4}
    Suppose $(f_1)-(f_4)$ hold. Then there exist a constants $c_0\in (0,c_b]$ and a sequence $\{u_n\} \subset \mathbb{X}_0^{s,2}(\Omega)$ satisfying 
    \begin{align}
        \Phi_b(u_n) \rightarrow c_0, ~\text{  and   } \left(1+\|u_n\|_{\mathbb{X}_0^{s,2}(\Omega)}\right)\|\Phi'_b(u_n)\|_{\mathbb{X}_0^{s,2}(\Omega)} \rightarrow 0.
    \end{align}
\end{lemma}
\begin{proof}
    From $(f_1)-(f_3)$ and \eqref{energy}, there exist $\mu_0>0$ and $\rho_0>0$ such that
    \begin{align}
        \Phi_b(u)\geq\rho_0,~\forall \,u\in \partial B_{\mu_0}(0),
    \end{align}
     where $B_{\mu_0}(0)=\{u\in{\mathbb{X}_0^{s,2}(\Omega)}:||u||_{\mathbb{X}_0^{s,2}(\Omega)}\leq \mu_0\}$. Let us consider $w_k\in \mathcal{N}_b$ such that
     \begin{align}\label{eq4.6}
         c_b\leq \Phi_b(w_k)< c_b+\frac{1}{k},~\forall~k\in \mathbb{N}.
     \end{align}
     Since $\Phi_b(tw_k)<0$ for large $t>0$, we have by the mountain pass theorem that there exists a sequence $\{u_{k,n}\}_{n\in \mathbb{N}}\subset \mathbb{X}_0^{s,2}(\Omega) $ such that
     \begin{align}\label{eq4.7}
         \Phi_b(u_{k,n}) \rightarrow c_k, ~\text{  and   } \|\Phi'_b(u_{k,n})\|_{\mathbb{X}_0^{s,2}(\Omega)}\left(1+\|u_{k,n}\|_{\mathbb{X}_0^{s,2}(\Omega)}\right) \rightarrow 0,
     \end{align}
     where $c_k \in[\rho_0, \sup_{t\geq 0}\Phi_b(tw_k)]$. Using Corollary \ref{cor2}, we get
     \begin{align}\label{eq4.8}
         \Phi_b(w_k)=\sup_{t\geq 0}\Phi_b(tw_k).
     \end{align}
     Thus, applying \eqref{eq4.6}, \eqref{eq4.7} and \eqref{eq4.8}, we obtain
     \begin{align}
         \Phi_b(u_{k,n}) \rightarrow c_k<c_b+\frac{1}{k}, ~\text{  and   } \|\Phi'_b(u_{k,n})\|_{\mathbb{X}_0^{s,2}(\Omega)}\left(1+\|u_{k,n}\|_{\mathbb{X}_0^{s,2}(\Omega)}\right) \rightarrow 0.
     \end{align}
     We now choose a sequence $\{n_k\}\subset \mathbb{N}$ such that
     \begin{align}
         \Phi_b(u_{k,n_k}) <c_b+\frac{1}{k}, ~\text{  and   } \|\Phi'_b(u_{k,n_k})\|_{\mathbb{X}_0^{s,2}(\Omega)}\left(1+\|u_{k,n_k}\|_{\mathbb{X}_0^{s,2}(\Omega)}\right) <\frac{1}{k}.
     \end{align}
     Let $u_k=u_{k,n_k}$ for $k\in \mathbb{N}$. Therefore, upto a subsequence, we have
     \begin{align}
          \Phi_b(u_{n}) \rightarrow c_0\in[\rho_0,c_b], ~\text{  and   } \|\Phi'_b(u_{n})\|_{\mathbb{X}_0^{s,2}(\Omega)}\left(1+\|u_{n}\|_{\mathbb{X}_0^{s,2}(\Omega)}\right) \rightarrow 0.
     \end{align}
     This completes the proof.
\end{proof} 

\section*{{\bf\emph{Proof of the Theorem \ref{T2.5}}}}
\begin{proof}
We divide the proof into two parts. The first part guarantees the existence of ground state energy for \eqref{New P} and in the next part, we prove $m_b\geq 2c_b$.\\

\noindent From Lemma \ref{lmn4}, there exists a sequence $\{u_n\} \subset \mathbb{X}_0^{s,2}(\Omega)$ such that
    \begin{align}\label{eq4.12}
        \Phi_b(u_n) \rightarrow c_0, ~\text{  and   } \langle \Phi'_b(u_n),u_n\rangle \rightarrow 0. 
    \end{align}
    Using \eqref{energy}, \eqref{eq21} and \eqref{eq4.12}, we obtain 
    \begin{align}
        c_0+1 \geq \Phi_b(u_n) -\frac{1}{4}\langle \Phi'_b(u_n),u_n \rangle \geq \frac{a(1-\theta)}{4} \|u_n\|_{\mathbb{X}_0^{s,2}(\Omega)}^2
    \end{align}
    for large $n\in \mathbb{N}$. Thus, we have $\{u_n\}$ is bounded in $\mathbb{X}_0^{s,2}(\Omega)$. Proceeding with a similar argument as in \cite{GCF2025} (see also \cite{T2014}), we deduce that there exists a $u\in \mathbb{X}_0^{s,2}(\Omega) \setminus \{0\}$ such that $\Phi'_b(u)=0$. This implies that $u\in \mathbb{N}_b$ is a nontrivial solution to the problem \eqref{New P} and $c_b \leq \Phi_b(u)$. On the other side, using \eqref{energy}, \eqref{eq21}, weak semicontinuity of norm and Fatou's lemma, we get
    \begin{align}
        c_b\geq c_0 =&\lim_{n\rightarrow\infty} \bigg[\Phi_b(u_n)-\frac{1}{4}\langle\Phi'_b(u_n),u_n\rangle\bigg]\nonumber \\
        =&\lim_{n\rightarrow\infty} \bigg[\frac{a}{4}\|u_n\|_{\mathbb{X}_0^{s,2}(\Omega)}^2+\int_\Omega \left(\frac{1}{4}f(u_n)u_n-F(u_n) \right) \d x \bigg]\nonumber \\
        \geq &\frac{a}{4} \liminf_{n\rightarrow\infty} \left(\|u_n\|_{\mathbb{X}_0^{s,2}(\Omega)}^2-\theta\lambda_1\|u_n\|^2_{L^2(\Omega)} \right)\nonumber \\
        &+\liminf_{n\rightarrow\infty}\int_\Omega \left(\frac{1}{4}f(u_n)u_n-F(u_n) +\frac{a\theta\lambda_1}{4}|u_n|^2\right) \d x \nonumber \\
        \geq &\frac{a}{4}  \left(\|u\|_{\mathbb{X}_0^{s,2}(\Omega)}^2-\theta\lambda_1\|u\|^2_{L^2(\Omega)} \right)+\int_\Omega \left(\frac{1}{4}f(u)u-F(u) +\frac{a\theta\lambda_1}{4}|u|^2\right) \d x \nonumber \\
         = &\frac{a}{4}  \|u\|_{\mathbb{X}_0^{s,2}(\Omega)}^2+\int_\Omega \left(\frac{1}{4}f(u)u-F(u) \right) \d x=\Phi_b(u)-\frac{1}{4}\langle\Phi'_b(u),u\rangle \nonumber \\
         =& \Phi_b(u) \geq \inf_{v\in \mathcal{N}_b}\Phi_b(v)=c_b.
    \end{align}
    Therefore, we have $\Phi_b(u)=c_b>0.$\\
    
 Now, from Theorem \ref{T2.4}, there exists a $u_b\in \mathcal{M}_b$ such that $\Phi_b(u_b)=m_b>0$. Moreover, using \eqref{energy}, Corollary \ref{cor} and Lemma \ref{lmn3}, we get
    \begin{align}
        m_b=&\Phi_b(u_b)=\sup_{k,l \geq 0}\Phi_b(ku_b^++lu_b^-) \nonumber \\
        =& \sup_{k,l\geq 0} \bigg\{\frac{a}{2}\left(k^2\|u_b^+\|_{\mathbb{X}_0^{s,2}(\Omega)}^2+l^2\|u_b^-\|_{\mathbb{X}_0^{s,2}(\Omega)}^2+2klA(u_b)\right) \nonumber \\
        &+\frac{b}{4}\left(k^2\|u_b^+\|_{\mathbb{X}_0^{s,2}(\Omega)}^2+l^2\|u_b^-\|_{\mathbb{X}_0^{s,2}(\Omega)}^2+2klA(u_b)\right)^2-\int_\Omega F(ku_b^+) \d x-\int_\Omega F(lu_b^-) \d x\bigg\} \nonumber \\
        > & \sup_{k,l \geq 0} \bigg\{\frac{a}{2}(k^2\|u_b^+\|_{\mathbb{X}_0^{s,2}(\Omega)}^2+\frac{a}{2}l^2\|u_b^-\|_{\mathbb{X}_0^{s,2}(\Omega)}^2 +\frac{b}{4}k^4\|u_b^+\|_{\mathbb{X}_0^{s,2}(\Omega)}^4+\frac{b}{4}l^4\|u_b^-\|_{\mathbb{X}_0^{s,2}(\Omega)}^4  \nonumber \\
        &-\int_\Omega F(ku_b^+) \d x-\int_\Omega F(lu_b^-) \d x \bigg\} \nonumber \\
        = & \sup_{k \geq 0}\Phi_b(ku^+_b)   + \sup_{l \geq 0} \Phi_b(lu_b^-) \geq 2c_b.
    \end{align}
    This completes the proof.
\end{proof}
\begin{remark}
    It is noteworthy to mention here that Theorem \ref{T2.5} can be proved by using the Nehari manifold technique combined with Lemma \ref{lmn5.6} and Lemma \ref{lmn5.7}.
 \end{remark}

\section{Asymptotic behavior of solutions}\label{sec5} 
In this Section, we prove the existence of least energy sign-changing solutions to the problem \eqref{New P2} and establish the convergence property of $u_b$ as $b \searrow 0$.

\section*{{\bf\emph{Proof of the Theorem \ref{T2.6}}}}
\begin{proof}
Note that the existence of sign-changing solutions to \eqref{New P2} can be obtained following similar arguments as in Section \ref{sec3} with $b=0$. Therefore, under the assumptions of Theorem \ref{T2.6}, there exists $v_0\in \mathcal{M}_0$ such that
    \begin{align}
        \Phi_0(v_0)=m_0=\inf_{u\in \mathcal{M}_0}\Phi_0(u)~\text{ and }~\Phi'_0(v_0)=0.
    \end{align}
    Thus, we conclude that problem \eqref{New P2} has a sign-changing solution. In particular, there exists a least energy sign-changing solution to \eqref{New P2}. \\

    \noindent We now focus on the asymptotic behavior of solutions $u_b$ as $b \searrow 0$. Recall Theorem \ref{T2.4}. Let $u_b\in \mathcal{M}_b$ be a sign-changing solution to the problem \eqref{New P} for any $b>0$ such that $\Phi_b(u_b)=m_b$. Let $w_1\in C_c^\infty(\Omega)$ be such that $w_1^\pm\neq 0$. Now, using the conditions $(f_1)-(f_3)$, we deduce that there exist $\gamma_1>0$ and 
    $$\gamma_2\geq \max
    \left\{ \frac{\|w_1^+\|^4_{\mathbb{X}_0^{s,2}(\Omega)}}{\|w_1^+\|_{4}^4},\frac{\|w_1^-\|^4_{\mathbb{X}_0^{s,2}(\Omega)}}{\|w_1^-\|_{4}^4}\right\},$$
    such that 
     \begin{align}\label{eq5.2}
        F(t)\geq \gamma_2 |t|^4-\gamma_1,~\forall ~t\in \mathbb{R}.
        \end{align}
    Therefore, for any $b\in [0,\kappa]$, where $\kappa$ is sufficiently small, using Lemma \ref{lmn5.5} combined with \eqref{eq5.2} and Young's inequality, we get
    \begin{align}\label{eq5.3}
        \Phi_b(u_b)=&m_b=\inf_{u\in {\mathbb{X}_0^{s,2}(\Omega)}, u^\pm \neq 0} \max_{k,l \geq 0}\Phi_b(ku^++lu^-) \nonumber\\
        \leq&  \max_{k,l \geq 0}\Phi_b(kw_1^++lw_1^-) \nonumber \\
        =& \max_{k,l\geq 0} \bigg\{\frac{a}{2}\left(k^2\|w_1^+\|_{\mathbb{X}_0^{s,2}(\Omega)}^2+l^2\|w_1^-\|_{\mathbb{X}_0^{s,2}(\Omega)}^2+2klA(w_1)\right)-\int_\Omega F(kw_1^+) \d x \nonumber \\
        &+\frac{b}{4}\left(k^2\|w_1^+\|_{\mathbb{X}_0^{s,2}(\Omega)}^2+l^2\|w_1^-\|_{\mathbb{X}_0^{s,2}(\Omega)}^2+2klA(w_1)\right)^2-\int_\Omega F(lw_1^-) \d x\bigg\}    \nonumber \\
        \leq &  \max_{k,l\geq 0} \bigg\{\frac{a}{2}k^2\|w_1^+\|_{\mathbb{X}_0^{s,2}(\Omega)}^2 +\frac{a}{2}l^2\|w_1^-\|_{\mathbb{X}_0^{s,2}(\Omega)}^2+aklA(w_1)+\frac{b}{4}k^4\|w_1^+\|_{\mathbb{X}_0^{s,2}(\Omega)}^4 \nonumber \\
        &+\frac{b}{4}l^4\|w_1^-\|_{\mathbb{X}_0^{s,2}(\Omega)}^4+bk^2l^2A(w_1)^2+\frac{b}{2}k^2l^2\|w_1^+\|_{\mathbb{X}_0^{s,2}(\Omega)}^2 \|w_1^-\|_{\mathbb{X}_0^{s,2}(\Omega)}^2\nonumber\\
        &+bk^3l\|w_1^+\|_{\mathbb{X}_0^{s,2}(\Omega)}^2 A(w_1)+bkl^3\|w_1^-\|_{\mathbb{X}_0^{s,2}(\Omega)}^2A(w_1) +2\gamma_1 |\Omega| \nonumber \\ 
        &-\gamma_2k^4 \|w_1^+\|_{4}^4 -\gamma_2l^4 \|w_1^-\|_{4}^4 \bigg\} \nonumber \\
\leq &  \max_{k,l\geq 0} \bigg\{\frac{a}{2}k^2\|w_1^+\|_{\mathbb{X}_0^{s,2}(\Omega)}^2 +\frac{a}{2}l^2\|w_1^-\|_{\mathbb{X}_0^{s,2}(\Omega)}^2+a\frac{(k^2+l^2)}{2}A(w_1)+\frac{b}{2}k^4\|w_1^+\|_{\mathbb{X}_0^{s,2}(\Omega)}^4 \nonumber \\
         & +\frac{b}{2}l^4\|w_1^-\|_{\mathbb{X}_0^{s,2}(\Omega)}^4 + \frac{b(k^4 A(w_1)^2+l^4 A(w_1)^2)}{2}+b\left(\frac{3k^4}{4}+\frac{l^4}{4}\right)\|w_1^+\|_{\mathbb{X}_0^{s,2}(\Omega)}^2 A(w_1) \nonumber \\
         & + b\left(\frac{k^4}{4}+\frac{3l^4}{4}\right) \|w_1^-\|_{\mathbb{X}_0^{s,2}(\Omega)}^2 A(w_1) +2\gamma_1|\Omega|-k^4\|w_1^+\|_{\mathbb{X}_0^{s,2}(\Omega)}^4-l^4\|w_1^-\|_{\mathbb{X}_0^{s,2}(\Omega)}^4 \bigg\} \nonumber \\
\leq&  \max_{k,l\geq 0} \bigg\{k^2\left(\frac{a}{2}\|w_1^+\|_{\mathbb{X}_0^{s,2}(\Omega)}^2+ \frac{a}{2}A(w_1)\right)-k^4\bigg(\|w_1^+\|_{\mathbb{X}_0^{s,2}(\Omega)}^4-\frac{\kappa}{2}\|w_1^+\|_{\mathbb{X}_0^{s,2}(\Omega)}^4\nonumber \\
     &-\frac{\kappa}{2}A(w_1)^2 -\frac{3\kappa}{4}\|w_1^+\|_{\mathbb{X}_0^{s,2}(\Omega)}^2 A(w_1) -\frac{\kappa}{4}\|w_1^-\|_{\mathbb{X}_0^{s,2}(\Omega)}^2 A(w_1)\bigg) +2\gamma_1 |\Omega| \nonumber \\
     &+l^2\left(\frac{a}{2}\|w_1^-\|_{\mathbb{X}_0^{s,2}(\Omega)}^2+ \frac{a}{2}A(w_1)\right)-l^4\bigg(\|w_1^-\|_{\mathbb{X}_0^{s,2}(\Omega)}^4-\frac{\kappa}{2}\|w_1^-\|_{\mathbb{X}_0^{s,2}(\Omega)}^4\nonumber \\
     &-\frac{\kappa}{2}A(w_1)^2 -\frac{\kappa}{4}\|w_1^+\|_{\mathbb{X}_0^{s,2}(\Omega)}^2 A(w_1) -\frac{3\kappa}{4}\|w_1^-\|_{\mathbb{X}_0^{s,2}(\Omega)}^2 A(w_1)\bigg) \bigg \}\nonumber\\
     &:= \aleph \in (0, \infty) .
    \end{align}
Let $\{b_n\}$ be a sequence such that $b_n\searrow 0$ as $n\rightarrow \infty$ and let $u_{b_n}$ be the corresponding weak solution to \eqref{New P}. Therefore, using \eqref{energy}, \eqref{eq21}, \eqref{eq5.3} and for sufficiently large $n\in \mathbb{N}$, we obtain
\begin{align}
    \aleph +1 \geq \Phi_{b_n}(u_{b_n})=\Phi_{b_n}(u_{b_n})-\frac{1}{4}\langle \Phi'_{b_n}(u_{b_n}),u_{b_n} \rangle \geq \frac{a(1-\theta)}{4} \|u_{b_n}\|_{\mathbb{X}_0^{s,2}(\Omega)}^2.
\end{align}
Thus, the sequence $\{u_{b_n}\}$ is bounded in $\mathbb{X}_0^{s,2}(\Omega)$. Therefore, there exists a subsequence $\{b_n\}$ (still denoted by $\{b_n\}$) and $u_0\in \mathbb{X}_0^{s,2}(\Omega)$ such that 
        \begin{align}
            u_{b_n}^\pm &\rightharpoonup u_0^\pm\in \mathbb{X}_0^{s,2}(\Omega)~\text{weakly in}~ \mathbb{X}_0^{s,2}(\Omega),\label{weak}\\
        u_{b_n}^\pm &\rightarrow u_0^\pm~\text{strongly in}~ L^r(\Omega)~\text{for}~r\in[1,2^*), \\ 
        u_{b_n}(x)& \rightarrow u_0(x)~a.e.~\text{in}~\Omega.
        \end{align}
        Proceeding with a similar idea from Lemma \ref{lmn5.6} and using \cite{S2015}, we deduce that
        \begin{align}
            u_{b_n}^\pm \rightarrow u_0^\pm~\text{strongly in}~ \mathbb{X}_0^{s,2}(\Omega)~\text{with}~u_0^\pm\neq0.
        \end{align}
        From \eqref{weak}, we get
        \begin{align}\label{5.10}
            \langle \Phi_0'(u_0),v \rangle=&a \langle u_0, v\rangle_{\mathbb{X}_0^{s,2}(\Omega)}- \int_\Omega f(u_0) v \d x \nonumber \\
            =& \lim_{n\rightarrow \infty}\bigg\{ \left(a+b_n\|u_{b_n}\|_{{\mathbb{X}_0^{s,2}(\Omega)}}\right)\langle u_{b_n}, v\rangle_{\mathbb{X}_0^{s,2}(\Omega)}- \int_\Omega f(u_{b_n}) v \d x  \bigg\} \nonumber \\
            =& \lim_{n\rightarrow \infty} \langle \Phi_{b_n}'(u_{b_n}),v \rangle=0,~\forall\, v\in C_c^\infty(\Omega).
        \end{align}
        Thus, we conclude that $\Phi_0'(u_0)=0$ and $u_0\in \mathcal{M}_0$. Therefore, we have $\Phi_0(u_0)\geq m_0.$ It remains to prove that $\Phi_0(u_0)\leq m_0.$ Recall that $v_0$ is a sign-changing solution to the problem \eqref{New P2}. Let $b_n \in [0,\kappa]$. Using $(f_3)$, there exists $j>0$ such that
        \begin{align}
            &\Phi_{b_n}(kv_0^++lv_0^-) \nonumber \\
            =&\frac{a}{2}\left(k^2\|v_0^+\|_{\mathbb{X}_0^{s,2}(\Omega)}^2+l^2\|v_0^-\|_{\mathbb{X}_0^{s,2}(\Omega)}^2+2klA(v_0)\right) \nonumber \\
        &+\frac{b_n}{4}\left(k^2\|v_0^+\|_{\mathbb{X}_0^{s,2}(\Omega)}^2+l^2\|v_0^-\|_{\mathbb{X}_0^{s,2}(\Omega)}^2+2klA(v_0)\right)^2-\int_\Omega F(kv_0^+) \d x-\int_\Omega F(lv_0^-) \d x  \nonumber \\
        \leq&  k^2\left(\frac{a}{2}\|v_0^+\|_{\mathbb{X}_0^{s,2}(\Omega)}^2+ \frac{a}{2}A(v_0)\right)+k^4\bigg(\frac{\kappa}{2}\|v_0^+\|_{\mathbb{X}_0^{s,2}(\Omega)}^4+\frac{\kappa}{2}A(v_0)^2 +\frac{3\kappa}{4}\|v_0^+\|_{\mathbb{X}_0^{s,2}(\Omega)}^2 A(v_0) \nonumber \\
     & +\frac{\kappa}{4}\|v_0^-\|_{\mathbb{X}_0^{s,2}(\Omega)}^2 A(v_0)\bigg) -\int_\Omega F(kv_0^+) \d x \nonumber \\
     &+l^2\left(\frac{a}{2}\|v_0^-\|_{\mathbb{X}_0^{s,2}(\Omega)}^2+ \frac{a}{2}A(v_0)\right)+l^4\bigg(\frac{\kappa}{2}\|v_0^-\|_{\mathbb{X}_0^{s,2}(\Omega)}^4 +\frac{\kappa}{2}A(v_0)^2 +\frac{\kappa}{4}\|v_0^+\|_{\mathbb{X}_0^{s,2}(\Omega)}^2 A(v_0)\nonumber \\
     & +\frac{3\kappa}{4}\|v_0^-\|_{\mathbb{X}_0^{s,2}(\Omega)}^2 A(v_0)\bigg) -\int_\Omega F(lv_0^-) \d x
        <0, ~\forall~ k+l \geq j.
        \end{align}
        Thus, from Lemma \ref{lmn5.4}, there exist $k_n, l_n>0$ such that $k_nv_0^++l_nv_0^-\in \mathcal{M}_{b_{n}}$. Using Corollary \eqref{cor}, Lemma \ref{lmn5.6} and \eqref{5.10}, we obtain $0<k_n,l_n<j.$ Now, using $\Phi_0'(v_0)=0$, \eqref{energy}, \eqref{energy 0} and Lemma \ref{lm 2.1}, we get
        \begin{align}\label{eq5.11}
            m_0=&\Phi_0(v_0)=\Phi_{b_{n}}(v_0) -\frac{b_n}{4} \|v_0\|^4_{\mathbb{X}_0^{s,2}(\Omega)} \nonumber \\
            \geq &\Phi_{b_n}(k_nv_0^++l_nv_0^-)+\frac{1-k_n^4}{4}\langle \Phi_{b_n}'(v_0),v_0^+\rangle+\frac{1-l_n^4}{4}\langle \Phi_{b_n}'(v_0),v_0^-\rangle -\frac{b_n}{4} \|v_0\|^4_{\mathbb{X}_0^{s,2}(\Omega)} \nonumber \\
            \geq & m_{b_n} -\frac{1+S^4}{4}|\langle \Phi_{b_n}'(v_0),v_0^+\rangle|-\frac{1+S^4}{4}|\langle \Phi_{b_n}'(v_0),v_0^-\rangle| -\frac{b_n}{4} \|v_0\|^4_{\mathbb{X}_0^{s,2}(\Omega)} \nonumber \\
            = & m_{b_n} -b_n\frac{1+S^4}{4} \|v_0\|_{\mathbb{X}_0^{s,2}(\Omega)}^2 \left(\|v_0^+\|_{\mathbb{X}_0^{s,2}(\Omega)}^2+A(v_0) \right) \nonumber \\
            &-b_n\frac{1+S^4}{4} \|v_0\|_{\mathbb{X}_0^{s,2}(\Omega)}^2 \left(\|v_0^-\|_{\mathbb{X}_0^{s,2}(\Omega)}^2+A(v_0) \right) -\frac{b_n}{4} \|v_0\|^4_{\mathbb{X}_0^{s,2}(\Omega)}.
        \end{align}
        Therefore, from \eqref{eq5.11}, we get
        \begin{align}\label{eq5.12}
            \limsup_{n\rightarrow \infty}m_{b_n} \leq m_0.
        \end{align}
        Therefore, using \eqref{energy}, \eqref{energy 0} and \eqref{eq5.12}, we conclude that
        \begin{align}
            m_0 \leq \Phi_0(u_0)=\limsup_{n\rightarrow \infty} \Phi_{b_n}(u_{b_n})=\limsup_{n\rightarrow \infty} m_{b_n} \leq m_0.
        \end{align}
        Hence, we have $\Phi_0(u_0)=m_0$. This completes the proof.
\end{proof}

\noindent \textbf{Conflict of interest statement:} On behalf of the authors, the corresponding author states that there is no conflict of interest.
\newline
\textbf{Data availability statement:} Data sharing does not apply to this article as no datasets were generated or analysed during the current study.
\newline
\textbf{Use of AI:} No AI tools is used for preparation of the manuscript.

\section*{Acknowledgement}
\noindent 
SB would like to thank the Council of Scientific and Industrial Research (CSIR), Govt. of India, for the financial assistance to carry out this research work [grant no. 09/0874(17164)/2023-EMR-I]. SG gratefully acknowledges financial support for this research under the ARG-MATRICS grant No. ANRF/ARGM/2025/001570/MTR, Anusandhan National Research Foundation (ANRF), Government of India. 
% SB and SG acknowledge the research facilities available at the Department of Mathematics, NIT Calicut.

%\bibliographystyle{plain}
%\bibliography{svik}

\end{document}